\documentclass[final,hidelinks,onefignum,onetabnum]{siamart251216}
\usepackage{amssymb,mathtools,bm,graphicx,placeins,enumitem,microtype}
\allowdisplaybreaks
\setlist{nosep,leftmargin=2em}
\newsiamthm{assumption}{Assumption}
\newsiamremark{remark}{Remark}
\crefname{assumption}{Assumption}{Assumptions}
\Crefname{assumption}{Assumption}{Assumptions}
\crefname{appendix}{Appendix}{Appendices}
\Crefname{appendix}{Appendix}{Appendices}
\newcommand{\R}{\mathbb{R}}

\newcommand{\Pp}{\mathbb{P}}
\newcommand{\E}{\mathbb{E}}
\newcommand{\B}{\mathcal{B}}
\newcommand{\X}{\mathcal{X}}

\newcommand{\Z}{\mathcal{Z}}
\newcommand{\G}{\mathcal{G}}
\newcommand{\F}{\mathfrak f}
\newcommand{\Rzero}{\mathcal{R}_0}

\newcommand{\ad}{\mathrm{ad}}
\newcommand{\prev}{\mathrm{prev}}

\newcommand{\TV}{\mathrm{TV}}
\newcommand{\Hell}{\mathrm{Hell}}

\newcommand{\dd}{\,\mathrm{d}}
\newcommand{\norm}[1]{\left\lVert #1\right\rVert}
\newcommand{\abs}[1]{\left\lvert #1\right\rvert}
\newcommand{\ip}[2]{\left\langle #1,#2\right\rangle}

\newcommand{\NumTrueLambda}{0.0100}
\newcommand{\NumPriorThresholdProb}{0.546}
\newcommand{\NumPrevThresholdProb}{0.999}
\newcommand{\NumBothThresholdProb}{0.995}
\newcommand{\NumPrevGrowthLower}{0.944}
\newcommand{\NumPrevGrowthUpper}{1.000}
\newcommand{\NumPrevALow}{0.289}
\newcommand{\NumPrevAHigh}{0.444}
\newcommand{\NumBothALow}{0.309}
\newcommand{\NumBothAHigh}{0.345}

\newcommand{\NumNonlinearPrevSlope}{0.99}
\newcommand{\NumNonlinearJointSlope}{0.98}
\newcommand{\NumSecondOrderSlope}{1.97}
\newcommand{\NumFirstOrderSlope}{0.97}

\headers{Bayesian epidemic-threshold inference}{Duc-Lam Duong}
\hypersetup{pdftitle={Bayesian inverse problems and epidemic-threshold uncertainty for a heterogeneous SIS reaction-diffusion model},pdfauthor={Duc-Lam Duong}}
\author{Duc-Lam Duong\thanks{Research Unit of Mathematical Sciences, University of Oulu, Finland
(\email{duc-lam.duong@oulu.fi}).}}

\title{Bayesian inverse problems and epidemic-threshold uncertainty
for a heterogeneous SIS reaction--diffusion model}
\begin{document}
\maketitle
\begin{abstract}
We study Bayesian inference of persistence from early observations of a
heterogeneous susceptible--infected--susceptible (SIS) reaction--diffusion
model with standard incidence. We show that the leading prevalence dynamics
depend on the unknown coefficients only through the difference between the
transmission and recovery rates and the infected diffusivity. The same pair determines the principal eigenvalue whose sign governs persistence, so early prevalence can inform the threshold even when the individual rates remain uncertain.
We prove global weak well-posedness and quantitative parameter stability
for the SIS system and use the resulting observation-map estimates to
establish Bayesian well-posedness.
For an initial infection of size $\varepsilon$, we derive a uniform
second-order expansion that describes how the separate rates enter
prevalence. Incidence depends explicitly on transmission at leading order.
When the noise standard deviation is proportional to $\varepsilon$, we
prove $O(\varepsilon)$ convergence in Hellinger distance of the nonlinear
posteriors to their low-prevalence limits for both prevalence alone and
joint prevalence and incidence data.
We bound the principal eigenvalue using the growth rate of the infected
$L^2$ norm between two times. This rate gives a lower bound for the full
nonlinear system, and the upper bound separates spatial transients from
nonlinear saturation. Uniform spectral estimates and bounds on infection
size and measurement error give sufficient conditions for persistence
and extinction.
We also bound the error in posterior persistence probabilities under
forward and eigenvalue approximation, accounting for posterior mass near
the threshold. Numerical examples illustrate threshold inference and the predicted \(O(\varepsilon)\) convergence of the nonlinear posteriors to their low-prevalence limits in Hellinger distance.
\end{abstract}

\begin{keywords}
Bayesian inverse problems, SIS epidemic model, reaction--diffusion systems,
weak identifiability, epidemic threshold, principal eigenvalue
\end{keywords}
\begin{MSCcodes}
62F15, 35K57, 35R30, 92D30, 65N30
\end{MSCcodes}

\section{Introduction}
At the beginning of an outbreak, measurements of the infected population can be available before the transmission and recovery rates are known. One question at this stage is whether the infection can persist under the prevailing conditions. An early assessment can inform decisions about surveillance and intervention before the long-term course of the outbreak is clear. In a heterogeneous population, persistence depends on where transmission and recovery occur and how infected individuals move. We study what can be learned about the persistence threshold from noisy observations over a short time interval and how the answer depends on what is measured.

We consider the susceptible--infected--susceptible (SIS) reaction--diffusion
model of Allen, Bolker, Lou, and Nevai \cite{AllenBolkerLouNevai2008}, in which
the transmission and recovery rates vary in space and the susceptible and
infected populations can have different diffusivities. The basic
reproduction number $\Rzero$ governs the stability of the disease-free
equilibrium and the existence of an endemic equilibrium
\cite{AllenBolkerLouNevai2008,Peng2009}. Huang et al.~\cite{HuangHanLiu2010}
study related thresholds with Dirichlet boundary conditions. Peng and Yi
\cite{PengYi2013} examine how spatial variation in epidemic risk and
population movement shape the endemic equilibrium. The connection between
reproduction numbers and principal eigenvalues is developed for more general
reaction--diffusion epidemic models in \cite{WangZhao2012,MagalWebbWu2019}.
Related SIS models incorporate advection \cite{CuiLou2016} and a varying
total population \cite{LiPengWang2017}. Hu, Jin, and Wang \cite{HuJinWang2024}
use the threshold to study how the location and size of a transmission
control zone affect persistence.

PDE-constrained optimization has been used to recover transmission and
recovery fields
\cite{XiangLiu2015,CoronelHuancasSepulveda2019,CoronelFrizHessZegarra2021}
and to identify diffusivities \cite{CoronelHuancasHessTello2024}.
Bayesian studies of other epidemic models have used data to estimate
heterogeneous parameters and quantify uncertainty in predictions
\cite{JhaCaoOden2020,ChenWuGhattas2021}. In the present setting, a short
record of noisy observations can be compatible with several combinations
of rates and diffusivities, which can imply different threshold signs.
We use a Bayesian formulation to combine the observations with prior
information on the coefficients and to propagate the resulting
uncertainty to the posterior probability that $\Rzero>1$.

The link between persistence and the individual rates can be seen in the
equation for the infected density, see \eqref{eq:SIS} below. Let \(S\) and \(I\) denote the susceptible and infected densities, respectively, and set \(N=S+I\). Using \(\beta SI/N=\beta I-\beta I^2/N\), we can write the infected equation \eqref{eq:SIS-I} as
\begin{equation}\label{eq:intro-I-rewrite}
 \partial_t I=d_I\Delta I+(\beta-\gamma)I-\beta I^2/N.
\end{equation}
When $I$ is small and $S$ remains positive, the leading dynamics of $I$
depend only on $q=\beta-\gamma$ and $d_I$. The same pair determines the
principal eigenvalue $\lambda_*$ of $d_I\Delta+q$ under Neumann boundary
conditions. The sign of $\lambda_*$ determines whether $\Rzero$ exceeds
one, by the equivalence of Allen et al.\
\cite[Lemma~2.3(d)]{AllenBolkerLouNevai2008} recalled in
\cref{thm:threshold-equivalence}. Thus early prevalence observations can
inform the threshold even when the separate rates remain uncertain.
In \eqref{eq:intro-I-rewrite}, $\beta$ enters separately through the
quadratic term. At low prevalence, the incidence $\beta SI/(S+I)$ is
close to $\beta I$, so incidence observations can provide information
about transmission at leading order.

In nonspatial epidemic models, the relation between early growth and the
reproduction number depends on the distribution of generation intervals
\cite{WallingaLipsitch2007}. In a spatial model, movement of infected
individuals into a region can raise a local count even when $\Rzero<1$.
The growth rate of the spatial $L^2$ norm over a short interval can also
differ from $\lambda_*$ because the initial infection excites several
spatial modes. The quadratic term in
\eqref{eq:intro-I-rewrite} can further reduce growth as prevalence increases.
We therefore ask how large a change in the infected density over a short
interval must be to determine the threshold sign, accounting for these
effects and for measurement error.

Our contributions are the following.
\begin{enumerate}[label=(\roman*)]
\item We prove global weak well-posedness and quantitative stability for
the SIS system when the rate fields and diffusivities vary
(\cref{thm:forward-wp,thm:forward-stability}). We use the resulting
observation-map estimates to establish Bayesian well-posedness
(\cref{sec:bayes}) and bound the error in the posterior persistence probability
under forward and eigenvalue approximation. 
\item We derive a uniform second-order expansion of the infected density
for an initial infection of size $\varepsilon$ (\cref{thm:second-order}).
With noise scaled by $\varepsilon$, we identify the coefficient of $\varepsilon^2$ in the Kullback--Leibler divergence between prevalence data laws for parameters
with the same $(q,d_I)$ (\cref{cor:KL-low-prevalence}).
\item Under the same noise scaling, we prove $O(\varepsilon)$ convergence
in Hellinger distance to the limiting prevalence and joint
prevalence--incidence posteriors, uniformly for bounded rescaled data
(\cref{thm:low-prevalence-posterior,cor:joint-limit}). 
\item We bound the principal eigenvalue using the growth of the infected
$L^2$ norm between two times (\cref{thm:growth-threshold}). Positive growth
of this norm implies persistence for the nonlinear system. The upper bound
accounts for spatial transients and nonlinear saturation and gives a
sufficient condition for extinction using uniform spectral estimates
(\cref{prop:uniform-spectral}) and bounds on infection size and
measurement error. 
We also bound the error in the posterior persistence probability under forward and eigenvalue approximation, accounting for posterior mass near the threshold (\cref{thm:threshold-approx}).
\end{enumerate}

The numerical experiments in \cref{sec:numerics} illustrate how prevalence
can inform persistence and how the nonlinear posteriors approach their
low-prevalence limits.

We use a common prior on the original coefficients in the nonlinear and
limiting experiments. The function-space formulation of Bayesian inversion
\cite{CotterDashtiRobinsonStuart2009,Stuart2010,DashtiStuart2017} lets us
control the difference between these posteriors in terms of the forward
approximation error \cite{CotterDashtiStuart2010,Sprungk2020}.
In the limiting prevalence experiment, the likelihood reweights the
distribution of $(q,d_I)$. The conditional-prior principle for partially
identified models \cite{Poirier1998,Gustafson2015} then describes the
posterior uncertainty in the rates by combining the updated law of this
pair with the conditional prior on the rates (\cref{prop:conditional-prior}).

The paper is organised as follows. We introduce the model and its threshold
in \cref{sec:model} and establish the forward estimates in \cref{sec:forward}.
We formulate the Bayesian inverse problem and construct admissible priors
in \cref{sec:bayes}. In \cref{sec:low-prevalence} we derive the low-prevalence
expansions, posterior limits, and growth bounds. We study posterior
uncertainty in the threshold in \cref{sec:threshold} and discuss the
results with numerical illustrations in \cref{sec:numerics}.
\Cref{app:spectral} gives the spectral estimates and their extension to
weighted observations. \Cref{app:prior-approximation} gives a concrete
prior construction and conditions for approximation, and
\cref{app:implementation} describes the numerical methods.

\section{Model, parameters, and observations}\label{sec:model}

Let $\Omega\subset\R^d$, $1\le d\le3$, be bounded and connected with $C^2$ boundary,
let $T>0$, and write $Q_T=\Omega\times(0,T)$. The susceptible and infected
densities $S$ and $I$ solve the system of Allen et al.~\cite{AllenBolkerLouNevai2008},
given by
\begin{subequations}\label{eq:SIS}
\begin{align}
  \partial_t S-d_S\Delta S
  &=-\beta(x)\F(S,I)+\gamma(x)I,
       &&\text{in }Q_T,\label{eq:SIS-S}\\
  \partial_t I-d_I\Delta I
  &= \beta(x)\F(S,I)-\gamma(x)I,
       &&\text{in }Q_T,\label{eq:SIS-I}\\
  \partial_\nu S=\partial_\nu I&=0,
       &&\text{on }\partial\Omega\times(0,T),\label{eq:SIS-bc}\\
  S(0)=S_0,\qquad I(0)&=I_0,
       &&\text{in }\Omega,\label{eq:SIS-ic}
\end{align}
\end{subequations}
where $\partial_\nu$ is the outward normal derivative and, for $s,i\ge0$,
the standard-incidence function is
\begin{equation}\label{eq:f-def}
  \F(s,i)=
  \begin{cases}
    \dfrac{si}{s+i},&s+i>0,\\[1ex]
    0,&s=i=0.
  \end{cases}
\end{equation}
Under local random mixing, the infection rate per susceptible is $\beta I/(S+I)$, proportional to the local fraction of infected individuals \cite[Section~2.1]{hethcote2000mathematics}. We refer to \eqref{eq:SIS-S} and \eqref{eq:SIS-I} as the susceptible and infected equations, respectively.

For fix positive constants $0<\underline\beta<\overline\beta, \, 0<\underline\gamma<\overline\gamma$ and $0<\underline d<\overline d$, 
we denote the ambient parameter space and admissible set as
\begin{align}
  \X&=C(\overline\Omega)\times C(\overline\Omega)\times\R^2,
     \label{eq:X-def}\\
  \X_\ad&=\bigl\{\theta=(\beta,\gamma,d_S,d_I)\in\X:
      \underline\beta\le\beta\le\overline\beta,
      \ \underline\gamma\le\gamma\le\overline\gamma,
      \ \underline d\le d_S,d_I\le\overline d\bigr\}.
      \label{eq:Xad-def}
\end{align}
Here and below, inequalities between functions are pointwise.  We equip $\X$ with
\begin{equation}\label{eq:X-norm}
  \norm{\theta}_{\X}
  =\norm{\beta}_{L^\infty}
   +\norm{\gamma}_{L^\infty}+|d_S|+|d_I|.
\end{equation}
With this norm, $\X$ is a separable Banach space, and the bounds in
\eqref{eq:Xad-def} define a closed bounded subset $\X_\ad \subset \X$.

For $\theta\in\X_\ad$ and $z=(S,I)$, write
\[
 \mathcal F_\theta(z)=\bigl(-\beta\F(S,I)+\gamma I,\,
                              \beta\F(S,I)-\gamma I\bigr).
\]
We call $\mathcal F_\theta$ the reaction map. Its components
$\mathcal F_{\theta,S}$ and $\mathcal F_{\theta,I}$ describe the local
changes in the densities due to infection and recovery.

We assume throughout that
\begin{equation}\label{eq:initial-assumption}
  S_0,I_0\in L^\infty(\Omega),
  \qquad S_0\ge s_*>0,
  \qquad I_0\ge0.
\end{equation}
The strict lower bound on $S_0$ controls the low-prevalence remainder.
Existence and uniqueness of nonnegative weak solutions hold for
nonnegative $L^\infty$ initial data.

The following elementary properties of $\F$ will be used repeatedly.

\begin{lemma}\label{lem:f-properties}
For $s,i\ge0$, we have
\begin{equation}\label{eq:f-bounds}
  0\le\F(s,i)\le\min\{s,i\}.
\end{equation}
The function $\F$ is globally Lipschitz on $[0,\infty)^2$ with constant
one in the Euclidean norm. On $(0,\infty)^2$, its partial derivatives
are given by
\begin{equation}\label{eq:f-derivatives}
  \partial_s\F(s,i)=\frac{i^2}{(s+i)^2},
  \qquad
  \partial_i\F(s,i)=\frac{s^2}{(s+i)^2}.
\end{equation}
Moreover, for every $m>0$, the first derivative of $\F$ is Lipschitz on
$\{(s,i):s,i\ge0,\ s+i\ge m\}$, with a constant depending only on $m$.
\end{lemma}

\begin{proof}
The bound \eqref{eq:f-bounds} follows from $s+i\ge\max\{s,i\}$.
Away from the origin,
$|\nabla\F(s,i)|^2=(s^4+i^4)/(s+i)^4\le1$.
Integrating along line segments gives the Lipschitz estimate, which
extends to the origin by continuity. Differentiating
\eqref{eq:f-derivatives} shows that the Hessian is bounded on
$\{s,i\ge0,\ s+i\ge m\}$, a convex set.
\end{proof}

Set
\begin{equation}\label{eq:Z-def}
  \Z_T=
  \left[C([0,T];L^2(\Omega))\cap L^2(0,T;H^1(\Omega))\right]^2
\end{equation}
with its natural norm. For $u,v\in H^1(\Omega)$, define the weak
Neumann Laplacian $\Delta_N u\in H^1(\Omega)^*$ by
\[
 \langle\Delta_N u,v\rangle=-\int_\Omega\nabla u\cdot\nabla v\dd x ,
\]
where $\langle\cdot,\cdot\rangle$ denotes duality between $H^1(\Omega)^*$
and $H^1(\Omega)$. This formulation incorporates the zero normal flux
by taking all $v\in H^1(\Omega)$ as test functions.
A weak solution of \eqref{eq:SIS} is a nonnegative $z=(S,I)\in\Z_T$
with $\partial_t z\in L^2(0,T;H^1(\Omega)^*)^2$ and the initial values
in \eqref{eq:SIS-ic}. It satisfies, for $k\in\{S,I\}$ and almost every $t$,
\[
 \langle\partial_t z_k(t),v\rangle
 +d_k\int_\Omega\nabla z_k(t)\cdot\nabla v\dd x
 =\int_\Omega \mathcal F_{\theta,k}(z(t))\,v\dd x ,
 \qquad v\in H^1(\Omega).
\]
Here $z_S=S$ and $z_I=I$.

To include prevalence and incidence data in one framework, define
\begin{equation}\label{eq:E-def}
  \mathcal E_T
  =C([0,T];L^2(\Omega))^3
\end{equation}
and let $\mathcal O:\mathcal E_T\to\R^J$ be bounded and linear.  If
$\mathcal S(\theta)=(S_\theta,I_\theta)$ denotes the solution of \eqref{eq:SIS}, the
parameter-to-observation map is
\begin{equation}\label{eq:G-def}
  \G(\theta)=\mathcal O\left(
  S_\theta,I_\theta,
  \beta\F(S_\theta,I_\theta)
  \right).
\end{equation}
This includes spatial averages of prevalence and incidence at selected times,
as well as space--time averages. The incidence component belongs to
$C([0,T];L^2)$ by the Lipschitz property of $\F$ and continuity of the states in
$L^2$. These bounded observations avoid spatial point evaluation. Throughout,
prevalence-type data refer to instantaneous counts or spatial averages of
the infected density $I$, while incidence is the rate of new infections.
The normalized fraction $I/(S+I)$ is a different observable.

\subsection{The epidemic threshold}\label{sec:threshold-definitions}
The spatial basic reproduction number of \eqref{eq:SIS} is the generalized Rayleigh
quotient
\begin{equation}\label{eq:R0-def}
 \Rzero(\theta)
 =\sup_{0\ne\varphi\in H^1(\Omega)}
 \frac{\displaystyle\int_\Omega\beta\varphi^2\dd x}
 {\displaystyle d_I\int_\Omega|\nabla\varphi|^2\dd x
  +\int_\Omega\gamma\varphi^2\dd x}.
\end{equation}
For the Allen--Bolker--Lou--Nevai model, this threshold governs disease-free
stability and the existence of endemic equilibria \cite{AllenBolkerLouNevai2008}.
We call $\{\Rzero>1\}$ the persistence event.  Its dynamical interpretation
presumes a nonzero initial infection: if $I_0\equiv0$, then $I(t)\equiv0$
regardless of the value of $\Rzero$.  A general
operator-theoretic treatment of reproduction numbers for reaction--diffusion epidemic
systems can be found in \cite{WangZhao2012,MagalWebbWu2019}.

The equation for a small infection in \eqref{eq:intro-I-rewrite} is
governed to leading order by $d_I\Delta+(\beta-\gamma)$ with Neumann
boundary conditions. Its largest eigenvalue gives the exponential growth
rate of the principal spatial mode. We denote it by
\begin{equation}\label{eq:lambda-def}
 \lambda_*(\theta)
 =\sup_{0\ne\varphi\in H^1(\Omega)}
 \frac{\displaystyle
       \int_\Omega(\beta-\gamma)\varphi^2\dd x
       -d_I\int_\Omega|\nabla\varphi|^2\dd x}
      {\displaystyle\int_\Omega\varphi^2\dd x}.
\end{equation}
Allen et al.~\cite[Section~2]{AllenBolkerLouNevai2008} write the time
dependence of the mode as $e^{-\lambda t}$, so their eigenvalue is
$-\lambda_*$. We recall their threshold equivalence in our sign
convention and record the properties of the principal eigenfunctions
used below.

\begin{theorem}\label{thm:threshold-equivalence}
For every $\theta\in\X_\ad$, the suprema in \eqref{eq:R0-def} and
\eqref{eq:lambda-def} are attained.  The corresponding eigenfunctions may be chosen
strictly positive on $\overline\Omega$, and the principal eigenvalues are simple.  In
addition,
\begin{equation}\label{eq:threshold-equivalence}
 \Rzero(\theta)>1
 \quad\Longleftrightarrow\quad
 \lambda_*(\theta)>0,
\end{equation}
with analogous equivalences for equality and for the reverse inequalities.
\end{theorem}

\begin{proof}
For \eqref{eq:lambda-def}, choose a maximizing sequence $(v_n)$ with
$\|v_n\|_{L^2}=1$, and denote its quotient values by $\ell_n$.
The gradient bound
$d_I\|\nabla v_n\|_{L^2}^2\le\|\beta-\gamma\|_{L^\infty}-\ell_n$
gives boundedness in $H^1(\Omega)$. Along a subsequence, $v_n$ converges
weakly in $H^1$ and strongly in $L^2$. Its limit has unit $L^2$ norm.
The potential term converges by strong $L^2$ convergence, and the
gradient norm is weakly lower semicontinuous. Thus the quotient at the
limit is at least $\lim_n\ell_n=\lambda_*$.

For \eqref{eq:R0-def}, set
$a(u,v)=d_I\int_\Omega\nabla u\cdot\nabla v\dd x+\int_\Omega\gamma uv\dd x$
and $b(u,v)=\int_\Omega\beta uv\dd x$. Normalize a maximizing sequence
by $a(v_n,v_n)=1$. The positive lower bounds for $d_I$ and $\gamma$
bound its $H^1$ norm. A weak $H^1$ and strong $L^2$ subsequential limit
$v$ satisfies $a(v,v)\le1$ and $b(v,v)=\Rzero$. The quotient on
constant functions is positive, so $v\ne0$ and $b(v,v)/a(v,v)\ge\Rzero$.
Both suprema are therefore attained.

Taking absolute values preserves both quotients, so choose nonnegative
maximizers $\phi_\lambda$ and $\phi_R$. For a quotient of symmetric
bilinear forms, differentiating $A(\phi+hv,\phi+hv)/B(\phi+hv,\phi+hv)$
at $h=0$ gives $A(\phi,v)=\rho B(\phi,v)$, where $\rho$ is the value of
the quotient at $\phi$. Here the two stationarity conditions, for all
$v\in H^1(\Omega)$, are
\[
 \begin{aligned}
  -d_I\int_\Omega\nabla\phi_\lambda\cdot\nabla v\dd x
     +\int_\Omega(\beta-\gamma)\phi_\lambda v\dd x
     &=\lambda_*\int_\Omega\phi_\lambda v\dd x,\\
  b(\phi_R,v)&=\Rzero a(\phi_R,v).
 \end{aligned}
\]
These are weak Neumann eigenvalue equations. Their weak Laplacians
belong to $L^2(\Omega)$, so Neumann regularity on the $C^2$ domain gives
$\phi_\lambda,\phi_R\in H^2(\Omega)$ \cite[Section~3.1.1]{Lunardi1995}.
The equations read
\begin{align}
 d_I\Delta\phi_\lambda+(\beta-\gamma)\phi_\lambda
 &=\lambda_*\phi_\lambda,\qquad \partial_\nu\phi_\lambda=0,
 \label{eq:lambda-eigenproblem}\\
 \Rzero(-d_I\Delta\phi_R+\gamma\phi_R)
 &=\beta\phi_R,\qquad \partial_\nu\phi_R=0.
 \label{eq:R0-eigenproblem}
\end{align}
The differential identities hold in $\Omega$ in $L^2$, and the normal
derivatives vanish on $\partial\Omega$ in the trace sense.

The lowest eigenvalue of $-d_I\Delta-(\beta-\gamma)$ is $-\lambda_*$.
The lowest eigenvalue of $-d_I\Delta+\gamma-\beta/\Rzero$ is zero, since
its quadratic form $a-b/\Rzero$ is nonnegative and vanishes at $\phi_R$.
For these scalar Neumann problems with bounded real potentials on a
connected domain, the lowest eigenvalue is simple and its eigenfunction
is strictly positive on $\overline\Omega$
\cite[Proposition~2.4 and Theorem~4.5]{ArendtTerElstGluck2020}.
The kernel of the second operator is precisely the eigenspace at $\Rzero$
of the generalized problem.

The numerator in \eqref{eq:lambda-def} is $b(v,v)-a(v,v)$. Thus
$\Rzero>1$ and $\lambda_*>0$ are both equivalent to the existence of
$v$ with $b(v,v)>a(v,v)$. If either supremum is at its threshold value,
$b(v,v)\le a(v,v)$ for every $v$, with equality at a nonzero maximizer.
This gives the equivalence at equality, and the remaining case gives
the equivalence for the strict reverse inequalities.
\end{proof}

\section{Analysis of the forward problem}\label{sec:forward}

\subsection{Global well-posedness and positivity}\label{sec:heat-evolution}

For a diffusivity $\kappa>0$, the Neumann realization $A_\kappa$ of
$\kappa\Delta$ on $L^2(\Omega)$ is
\[
 A_\kappa v=\kappa\Delta v,\qquad
 D(A_\kappa)=\{v\in H^2(\Omega):\partial_\nu v=0
                             \text{ on }\partial\Omega\}.
\]
It is the restriction of $\kappa\Delta_N$ to functions for which the
weak Laplacian is in $L^2(\Omega)$. The symmetric Neumann form gives
$-A_\kappa$ as a nonnegative self-adjoint operator. Its heat semigroup
$P_\kappa(t)=e^{tA_\kappa}$ is an analytic contraction semigroup on
$L^2(\Omega)$, and $P_\kappa(t)v_0$ is the homogeneous Neumann heat
evolution from $v_0$ \cite[Chapter~3]{Lunardi1995}. The parabolic
comparison principle implies that $P_\kappa(t)$ preserves order. Since
it preserves constants, comparison with $\pm\|v_0\|_{L^\infty}$ gives
$\|P_\kappa(t)v_0\|_{L^\infty}\le\|v_0\|_{L^\infty}$. Interpolation gives
contractivity on $L^p(\Omega)$ for $2\le p\le\infty$.
For $t>0$ and $v_0\ge0$, $v_0\not\equiv0$, the continuous representative
of $P_\kappa(t)v_0$ is strictly positive on $\overline\Omega$. We refer
to this property as positivity improvement
\cite[Theorem~3.1 and Section~4.2]{ArendtTerElstGluck2020}.

For the inhomogeneous heat equation with forcing $g$, the
variation-of-constants formula is
\[
 v(t)=P_\kappa(t)v(0)+\int_0^t P_\kappa(t-r)g(r)\dd r .
\]
This identity in $L^2(\Omega)$ defines the mild solution for
$g\in L^1(0,T;L^2(\Omega))$. For a bounded multiplication potential $V$,
the same formula holds with the semigroup $e^{t(A_\kappa+V)}$
\cite[Chapters~3 and~4]{Lunardi1995}.

For global classical solvability of SIS systems with standard incidence,
see \cite[Section~1]{AllenBolkerLouNevai2008} and
\cite[Theorem~2.1]{CoronelHuancasHessTello2024}. We work with bounded
initial densities and give bounds uniform over the admissible parameters.

\begin{theorem}\label{thm:forward-wp}
Let \eqref{eq:initial-assumption} hold and let $\theta\in\X_\ad$.  Then
\eqref{eq:SIS} has a unique nonnegative global weak solution
$(S_\theta,I_\theta)\in\Z_T$. It satisfies, almost everywhere in $Q_T$,
\begin{align}
  I_\theta(x,t)
  &\le e^{\overline\beta t}\norm{I_0}_{L^\infty},
     \label{eq:I-upper}\\
  S_\theta(x,t)
  &\ge s_*e^{-\overline\beta t},
     \label{eq:S-lower}\\
  S_\theta(x,t)
  &\le \norm{S_0}_{L^\infty}
  +\overline\gamma\norm{I_0}_{L^\infty}
    \int_0^t e^{\overline\beta r}\dd r.
     \label{eq:S-upper}
\end{align}
Furthermore,
\begin{equation}\label{eq:mass-conservation}
  \int_\Omega(S_\theta+I_\theta)(x,t)\dd x
  =\int_\Omega(S_0+I_0)(x)\dd x,
  \qquad 0\le t\le T.
\end{equation}
The constants in the $\Z_T$ and $L^\infty(Q_T)^2$ bounds are uniform over
$\theta\in\X_\ad$.

If, in addition, $0<\alpha<1$, $\partial\Omega$ is of class $C^{2+\alpha}$,
$\beta,\gamma\in C^\alpha(\overline\Omega)$,
$S_0,I_0\in C^{2+\alpha}(\overline\Omega)$, and the Neumann compatibility conditions
hold, then the solution is classical and belongs to
$C^{2+\alpha,1+\alpha/2}(\overline\Omega\times[0,T])^2$.
\end{theorem}

\begin{proof}
By \cref{lem:f-properties}, $\mathcal F_\theta$ is globally Lipschitz
on the nonnegative cone of $L^2(\Omega)^2$, the set of pairs with both
components nonnegative almost everywhere. Write $v^+=\max\{v,0\}$ and
define a reaction map on all of $L^2(\Omega)^2$ by
\[
 \widetilde{\mathcal F}_\theta(S,I)=\mathcal F_\theta(S^+,I^+).
\]
The positive-part map is $1$-Lipschitz in $L^2$, so the extended
reaction is globally Lipschitz. Let
$P(t)=\operatorname{diag}(P_{d_S}(t),P_{d_I}(t))$ and $z_0=(S_0,I_0)$.
Replacing the reaction in the two equations by
$\widetilde{\mathcal F}_\theta$, the mild formulation is
\[
 z(t)=P(t)z_0+\int_0^t P(t-r)\widetilde{\mathcal F}_\theta(z(r))\dd r .
\]
The integral map on $C([0,\tau];L^2(\Omega)^2)$ has Lipschitz constant
at most $\tau$ times that of $\widetilde{\mathcal F}_\theta$.
It is a contraction for small $\tau$, and the same step length can
be used from each new initial value. Iteration gives a mild solution
on any finite interval. Its forcing
$\widetilde{\mathcal F}_\theta(z)\in C([0,T];L^2(\Omega)^2)$ can be
treated as given in the linear heat equations. The linear energy
estimate gives $z\in\Z_T$ and
$\partial_t z\in L^2(0,T;H^1(\Omega)^*)^2$.
Indeed, one can project the linear equations on the Neumann
eigenfunctions, test each projected equation with the projected
solution, and pass to the limit in the weak identities of
\cref{sec:model}.

The reaction is quasipositive, meaning that each component is
nonnegative when its own density is zero and the other is nonnegative.
Here $\mathcal F_{\theta,S}(0,I)=\gamma I$ and
$\mathcal F_{\theta,I}(S,0)=0$. To use these signs in the weak
equations, set $v^-=(-v)^+$ and test the $k$th equation with
$-z_k^-(t)$. The chain rule and the gradient term give
\[
 \frac12\frac{\dd}{\dd t}\|z_k^-(t)\|_{L^2}^2
   +d_k\|\nabla z_k^-(t)\|_{L^2}^2
 =-\int_\Omega \widetilde{\mathcal F}_{\theta,k}(z(t))z_k^-(t)\dd x
 \le0 .
\]
On the support of $z_k^-$, the $k$th input of the extended reaction
is zero. The initial negative parts vanish, so both components stay
nonnegative and the extended reaction agrees with the original one.

Since $\F(S,I)\le I$, the infected equation implies
\[
  \partial_t I-d_I\Delta I\le\overline\beta I.
\]
Comparison with the spatially constant supersolution
$e^{\overline\beta t}\|I_0\|_{L^\infty}$ gives \eqref{eq:I-upper}. Likewise,
$\F(S,I)\le S$ yields
\[
  \partial_t S-d_S\Delta S
  \ge-\overline\beta S,
\]
and comparison with $s_*e^{-\overline\beta t}$ gives \eqref{eq:S-lower}.
Finally,
\[
  \partial_t S-d_S\Delta S\le\overline\gamma I,
\]
so the variation-of-constants formula and order preservation of
$P_{d_S}(t)$ give an upper bound with forcing $\overline\gamma I$.
Applying the $L^\infty$ contraction estimate and \eqref{eq:I-upper}
gives \eqref{eq:S-upper}.

The bounds \eqref{eq:I-upper}--\eqref{eq:S-upper} bound the reaction
uniformly in $L^2(Q_T)^2$. Using $S$ and $I$ as test functions and
integrating the energy identities gives a uniform bound in $\Z_T$.
For uniqueness, subtract the weak equations for two nonnegative
solutions and test with their difference. The Lipschitz estimate for
$\mathcal F_\theta$ and Gronwall's inequality give equality of the
solutions because the initial difference is zero.

The conservation law \eqref{eq:mass-conservation} follows by adding the two
variational equations with the constant test function one, since the
reaction terms cancel. Under the additional H\"older and compatibility
assumptions, linear parabolic regularity with bounded forcing first
gives $C^{\alpha,\alpha/2}$ regularity. By \cref{lem:f-properties},
$\F(S,I)$ has the same H\"older regularity as the states. Since the
coefficients are H\"older continuous in space, the reaction terms are
also in $C^{\alpha,\alpha/2}$. The Neumann Schauder estimates yield the stated
classical regularity \cite[Sections~5.1.2 and~7.1.3]{Lunardi1995}.
\end{proof}

\begin{remark}\label{rem:positive-S}
The explicit bound \eqref{eq:S-lower} implies
\begin{equation}\label{eq:N-lower}
  S_\theta+I_\theta\ge s_*e^{-\overline\beta T}=:n_T>0
  \quad\text{on }Q_T.
\end{equation}
This lower bound controls the denominator of the quadratic infection term
in the low-prevalence analysis of \cref{sec:low-prevalence}.
\end{remark}

\subsection{Quantitative stability with respect to the parameters}

For fixed unit diffusivities, \cite{CoronelHuancasSepulveda2019} gives
an $L^2$ estimate for state differences under perturbations of the
transmission and recovery rates. We allow both diffusivities to vary
and also control the spatial gradients in $L^2(Q_T)$.

\begin{theorem}\label{thm:forward-stability}
Let $\theta_j=(\beta_j,\gamma_j,d_{S,j},d_{I,j})\in\X_\ad$ and let
$z_j=(S_j,I_j)=\mathcal S(\theta_j)$, $j=1,2$.  There exists
$C=C(T,\Omega,\X_\ad,S_0,I_0)>0$ such that
\begin{equation}\label{eq:forward-stability}
  \norm{z_1-z_2}_{\Z_T}
  \le C\norm{\theta_1-\theta_2}_{\X}.
\end{equation}
If $d_{S,1}=d_{S,2}$ and $d_{I,1}=d_{I,2}$, then, for every
$2\le p<\infty$, there exists $C_p<\infty$ such that
\begin{equation}\label{eq:forward-stability-Lp}
 \norm{z_1-z_2}_{C([0,T];L^p)^2}
 \le C_p\left(
 \norm{\beta_1-\beta_2}_{L^\infty}
 +\norm{\gamma_1-\gamma_2}_{L^\infty}
 \right).
\end{equation}
With $C_\infty$ in place of $C_p$, the right-hand side of
\eqref{eq:forward-stability-Lp} also bounds $\norm{z_1-z_2}_{L^\infty(Q_T)^2}$.
\end{theorem}

\begin{proof}
Let $e=z_1-z_2$ and $g=\mathcal F_{\theta_1}(z_1)-\mathcal F_{\theta_2}(z_2)$.
For each component $k\in\{S,I\}$, subtract the weak identities of
\cref{sec:model}. For almost every $t$ and every $v\in H^1(\Omega)$ this gives
\[
 \begin{aligned}
 \langle\partial_t e_k,v\rangle+d_{k,1}\int_\Omega\nabla e_k\cdot\nabla v\dd x
 =\int_\Omega g_kv\dd x \quad-(d_{k,1}-d_{k,2})\int_\Omega\nabla z_{2,k}\cdot\nabla v\dd x .
 \end{aligned}
\]
The last term pairs $(d_{k,1}-d_{k,2})\Delta_N z_{2,k}$ with $v$.
By the Cauchy--Schwarz inequality its absolute value is at most
\[
 |d_{k,1}-d_{k,2}|\|\nabla z_{2,k}\|_{L^2}\|\nabla v\|_{L^2}.
\]
Thus $z_{2,k}(t)\in H^1(\Omega)$ suffices for this term to act in the dual.
Let $\delta_b=\|\beta_1-\beta_2\|_{L^\infty}+\|\gamma_1-\gamma_2\|_{L^\infty}$.
The Lipschitz estimate for the reaction and the uniform state bounds in
\cref{thm:forward-wp} give $|g|\le C(|e|+\delta_b)$ almost everywhere.

Use $v=e_k(t)$ in the $k$th identity, which is allowed for almost
every $t$ since $e\in L^2(0,T;H^1(\Omega))^2$. The time derivative
pairing is $\frac12\frac{\dd}{\dd t}\|e_k(t)\|_{L^2}^2$ by the energy
chain rule. Adding the two identities and using Young's inequality
for the diffusion differences yields
\begin{align*}
 \frac12\frac{\dd}{\dd t}
 \bigl(\norm{e_S}_{L^2}^2+\norm{e_I}_{L^2}^2\bigr)
 &+\frac{\underline d}{2}
 \bigl(\norm{\nabla e_S}_{L^2}^2+\norm{\nabla e_I}_{L^2}^2\bigr)\\
 &\le C\bigl(\norm{e_S}_{L^2}^2+\norm{e_I}_{L^2}^2\bigr)
   +C\delta_b^2\\
 &\quad+C|d_{S,1}-d_{S,2}|^2\norm{\nabla S_2}_{L^2}^2
       +C|d_{I,1}-d_{I,2}|^2\norm{\nabla I_2}_{L^2}^2.
\end{align*}
Since $e(0)=0$, Gronwall's inequality and the uniform energy bound for
$z_2$ bound the supremum of $\|e(t)\|_{L^2}^2$. Integrating the same
inequality bounds the time integral of $\|\nabla e\|_{L^2}^2$, giving
\eqref{eq:forward-stability}.

When the diffusivities agree, the weak difference equation has forcing
$g$ and zero initial data. The variation-of-constants formula and the
$L^p$ contraction estimates of \cref{sec:heat-evolution} give
\[
 \|e(t)\|_{L^p(\Omega)^2}
 \le\int_0^t \|g(r)\|_{L^p(\Omega)^2}\dd r
 \le C_p\int_0^t\bigl(\|e(r)\|_{L^p(\Omega)^2}+\delta_b\bigr)\dd r,
 \qquad 2\le p\le\infty.
\]
Gronwall's inequality gives the claimed bounds. For $p<\infty$,
continuity in $L^p$ follows by interpolating continuity in $L^2$ with
the uniform $L^\infty$ bound. For $p=\infty$ we use the essential
supremum over $Q_T$.
\end{proof}

\begin{proposition}\label{prop:G-regularity}
The map $\G$ in \eqref{eq:G-def} is bounded and Lipschitz on $\X_\ad$:
\begin{equation}\label{eq:G-bounds}
 |\G(\theta)|\le M_\G,\qquad
 |\G(\theta_1)-\G(\theta_2)|\le L_\G\|\theta_1-\theta_2\|_\X.
\end{equation}
\end{proposition}
\begin{proof}
The estimates follow from \cref{thm:forward-wp,thm:forward-stability} and
the boundedness of $\mathcal O$. For the incidence component, we also use
the Lipschitz property of $\F$ to bound $\beta\F(S,I)$ and its parameter
differences in $C([0,T];L^2)$.
\end{proof}

\section{The Bayesian inverse problem}\label{sec:bayes}

\subsection{Posterior definition and well-posedness}

Let $\B(\X_\ad)$ denote the Borel $\sigma$-algebra inherited from $\X$.  The
results in this section hold for any prior satisfying the following support
condition. Concrete constructions are given in \cref{sec:prior-modelling}.

\begin{assumption}\label[assumption]{ass:prior}
$\mu_0$ is an arbitrary Borel probability measure on $\X$ satisfying
$\mu_0(\X_\ad)=1$.
\end{assumption}

Let the data be generated by
\begin{equation}\label{eq:data-model}
  y=\G(\theta)+\eta,
  \qquad \eta\sim N(0,\Gamma),
\end{equation}
where $\Gamma\in\R^{J\times J}$ is symmetric and positive definite.  The negative
log-likelihood, up to an additive constant, is
\begin{equation}\label{eq:Phi}
  \Phi(\theta;y)
  =\frac12\abs{\Gamma^{-1/2}(y-\G(\theta))}^2.
\end{equation}
For probability measures $\mu$ and $\nu$ dominated by a common measure $\vartheta$, we use
the convention
\begin{equation}\label{eq:Hellinger-def}
 d_\Hell(\mu,\nu)^2
 =\frac12\int
 \left(
   \sqrt{\frac{\dd\mu}{\dd\vartheta}}
  -\sqrt{\frac{\dd\nu}{\dd\vartheta}}
 \right)^2\dd\vartheta.
\end{equation}

We also use the total-variation convention
\begin{equation}\label{eq:TV-def}
 d_\TV(\mu,\nu)=\sup_A|\mu(A)-\nu(A)|
 =\frac12\int\left|\frac{\dd\mu}{\dd\vartheta}
                       -\frac{\dd\nu}{\dd\vartheta}\right|\dd\vartheta,
\end{equation}
for any common dominating measure $\vartheta$.

\begin{theorem}\label{thm:posterior-existence}
Under \cref{ass:prior}, Bayes' formula defines, for every $y\in\R^J$, a
posterior probability measure $\mu^y$ on $\X_\ad$ uniquely specified by
\begin{equation}\label{eq:posterior-RN}
  \frac{\dd\mu^y}{\dd\mu_0}(\theta)
  =\frac1{Z(y)}e^{-\Phi(\theta;y)},
  \qquad
  Z(y)=\int_{\X_\ad}e^{-\Phi(\theta;y)}\dd\mu_0(\theta).
\end{equation}
Moreover, $0<Z(y)\le1$.
\end{theorem}

\begin{proof}
Continuity of $\G$ makes the likelihood measurable.
It is strictly positive and at most one, so $0<Z(y)\le1$.
Dividing by $Z(y)$ gives a density of total mass one relative to $\mu_0$.
For general results on Bayesian well-posedness, see \cite{Stuart2010,Latz2020}.
\end{proof}

Regular conditional distributions are unique up to sets of data of
prior-predictive measure zero. Formula \eqref{eq:posterior-RN} specifies a
version for every $y$. The next theorem bounds changes in this posterior
when the data or the parameter-to-observation map vary. We use the same
estimate for posterior approximation and for the low-prevalence limit.

\begin{theorem}\label{thm:posterior-stability}
Let $\G_1,\G_2:\X_\ad\to\R^J$ be Borel maps satisfying
\begin{equation}\label{eq:G-uniform-M}
  \sup_{\theta\in\X_\ad}\max\{|\G_1(\theta)|,|\G_2(\theta)|\}\le M.
\end{equation}
Let $\mu_j^{y_j}$ be the posterior obtained from $\G_j$, data $y_j$, covariance
$\Gamma$, and the same prior $\mu_0$.  For every $R<\infty$ there is
$C=C(R,M,\Gamma)<\infty$ such that, whenever $|y_1|,|y_2|\le R$,
\begin{equation}\label{eq:posterior-stability-general}
 d_\Hell(\mu_1^{y_1},\mu_2^{y_2})
 \le C\left(
 |y_1-y_2|
 +\norm{\G_1-\G_2}_{L^2(\mu_0;\R^J)}
 \right).
\end{equation}
In particular, for the SIS forward map $\G$,
\begin{equation}\label{eq:data-Hell-Lipschitz}
 d_\Hell(\mu^{y_1},\mu^{y_2})\le C_R|y_1-y_2|,
 \qquad |y_1|,|y_2|\le R.
\end{equation}
Thus the Bayesian inverse problem is locally Lipschitz well posed in the Hellinger
metric.
\end{theorem}

\begin{proof}
Write $\Phi_j=\frac12|\Gamma^{-1/2}(y_j-\G_j)|^2$ and
$Z_j=\int e^{-\Phi_j}\dd\mu_0$. Uniform boundedness gives
\begin{equation}\label{eq:Z-lower}
 e^{-K}\le Z_j\le1,\qquad K=\tfrac12\|\Gamma^{-1/2}\|^2(R+M)^2.
\end{equation}
The difference-of-squares identity yields
\begin{equation}\label{eq:Phi-difference}
 |\Phi_1-\Phi_2|\le C\bigl(|y_1-y_2|+|\G_1-\G_2|\bigr).
\end{equation}
For nonnegative $a,b$, $|e^{-a}-e^{-b}|\le|a-b|$ and
$|e^{-a/2}-e^{-b/2}|\le|a-b|/2$. Hence
$|Z_1-Z_2|\le\|\Phi_1-\Phi_2\|_{L^1(\mu_0)}$ and the $L^2(\mu_0)$
difference of the unnormalized square-root densities is at most
$\|\Phi_1-\Phi_2\|_{L^2(\mu_0)}/2$.
Following \cite[proof of Theorem~4.2]{Stuart2010}, we separate the
change in the likelihood from the change in its normalizing constant.
The triangle inequality gives
\[
 \begin{aligned}
 \sqrt{2}\,d_\Hell(\mu_1^{y_1},\mu_2^{y_2})
 &\le Z_1^{-1/2}\|e^{-\Phi_1/2}-e^{-\Phi_2/2}\|_{L^2(\mu_0)}\\
 &\quad+|Z_1^{-1/2}-Z_2^{-1/2}|\|e^{-\Phi_2/2}\|_{L^2(\mu_0)} .
 \end{aligned}
\]
By \eqref{eq:Z-lower}, both normalizing constants lie in $[e^{-K},1]$,
where $z\mapsto z^{-1/2}$ is Lipschitz. The preceding bounds give
\[
 d_\Hell(\mu_1^{y_1},\mu_2^{y_2})\le C\|\Phi_1-\Phi_2\|_{L^2(\mu_0)}.
\]
Applying \eqref{eq:Phi-difference} completes the proof.
\end{proof}

\subsection{Approximation of posterior measures}

Let $\mathfrak a$ denote a collective approximation index.  It may include a truncation of the
unknown coefficients, a spatial mesh size, and a time step.  Let
$\G_{\mathfrak a}:\X_\ad\to\R^J$ be a measurable approximate forward map and define
\begin{equation}\label{eq:posterior-approx}
 \frac{\dd\mu_{\mathfrak a}^y}{\dd\mu_0}(\theta)
 =\frac1{Z_{\mathfrak a}(y)}
 \exp\left(-\frac12|\Gamma^{-1/2}(y-\G_{\mathfrak a}(\theta))|^2\right).
\end{equation}

\begin{theorem}\label{thm:posterior-approximation}
Assume that
\begin{equation}\label{eq:Ga-bounded}
 \sup_{\mathfrak a}\sup_{\theta\in\X_\ad}|\G_{\mathfrak a}(\theta)|<\infty
\end{equation}
and that, for a deterministic sequence $\psi_{\mathfrak a}\to0$,
\begin{equation}\label{eq:Ga-error}
 \norm{\G-\G_{\mathfrak a}}_{L^2(\mu_0;\R^J)}\le C\psi_{\mathfrak a}.
\end{equation}
Then, for each $R<\infty$,
\begin{equation}\label{eq:posterior-approx-rate}
 \sup_{|y|\le R}d_\Hell(\mu^y,\mu_{\mathfrak a}^y)
 \le C_R\psi_{\mathfrak a}.
\end{equation}
If $Q:\X_\ad\to\R$ is bounded and Borel measurable, then
\begin{equation}\label{eq:QoI-expectation-Hell}
 \sup_{|y|\le R}
 \left|\int Q\dd\mu^y-\int Q\dd\mu_{\mathfrak a}^y\right|
 \le2\sqrt2\norm{Q}_{L^\infty}C_R\psi_{\mathfrak a}.
\end{equation}
\end{theorem}

\begin{proof}
For approximation by perturbing the likelihood, see
\cite[Theorem~4.6 and Corollary~4.9]{Stuart2010}. In the present setting,
\eqref{eq:posterior-approx-rate} follows from \cref{thm:posterior-stability}
with $y_1=y_2=y$ and maps $\G$ and $\G_{\mathfrak a}$.
For the second claim, recall that
\[
 d_\TV(\mu,\nu)\le\sqrt2d_\Hell(\mu,\nu)
\]
under convention \eqref{eq:Hellinger-def}, and
$|\int Q\dd\mu-\int Q\dd\nu|\le2\norm{Q}_\infty d_\TV(\mu,\nu)$.
\end{proof}

The abstract forward condition \eqref{eq:Ga-error} follows directly from a state-space
error estimate.

If a state solver and coefficient approximation give
\begin{equation}\label{eq:state-solver-error}
 \|z-z_{\mathfrak a}\|_{C([0,T];L^2)^2}
 +\|\beta\F(S,I)-\beta_{\mathfrak a}\F(S_{\mathfrak a},I_{\mathfrak a})\|_{C([0,T];L^2)}
 \le C\psi_{\mathfrak a}
\end{equation}
uniformly in the parameters, boundedness of $\mathcal O$ gives
\eqref{eq:Ga-error} and uniform boundedness of $\G_{\mathfrak a}$. The same conclusion
holds with an $L^2(\mu_0)$ state-error bound and a separately bounded approximate
map. Specific discretization rates require additional state regularity and
integrable error constants. For a sufficient condition yielding a
combined coefficient, space, and time approximation rate, we refer
to \cref{app:forward}.

\subsection{Admissible prior modelling}\label{sec:prior-modelling}
Let $\nu_0$ be a Gaussian Borel probability measure on
$\mathcal U=C(\overline\Omega)^2\times\R^2$. In particular, its field components
have continuous sample paths. For $a<b$ set
\begin{equation}\label{eq:bounded-link}
 \ell(t)=(1+e^{-t})^{-1},\qquad L_{a,b}(t)=a+(b-a)\ell(t).
\end{equation}
Applying these links componentwise defines
\begin{equation}\label{eq:prior-transform}
 \mathcal T(u)=\bigl(L_{\underline\beta,\overline\beta}(u_\beta),
 L_{\underline\gamma,\overline\gamma}(u_\gamma),
 L_{\underline d,\overline d}(u_S),L_{\underline d,\overline d}(u_I)\bigr).
\end{equation}
Known diffusivities can be fixed. Since $0<\ell'\le1/4$, $\mathcal T$ is Lipschitz and $\mu_0=\mathcal T_\sharp\nu_0$ satisfies
\cref{ass:prior}. The latent posterior with density proportional to
$e^{-\Phi(\mathcal T(u);y)}$ pushes forward to \eqref{eq:posterior-RN} by the
change-of-measure identity. Continuous Mat\'ern fields are one possible
choice of latent prior \cite{LindgrenRueLindstrom2011}.

Early prevalence suggests eliciting $q=\beta-\gamma$ and conditional recovery
uncertainty. To respect the coefficient bounds in \eqref{eq:Xad-def}, define
\begin{align}\label{eq:q-feasible}
 q_-&=\underline\beta-\overline\gamma,\quad
 q_+=\overline\beta-\underline\gamma,\notag\\
 \gamma_-(q)&=\max\{\underline\gamma,\underline\beta-q\},\quad
 \gamma_+(q)=\min\{\overline\gamma,\overline\beta-q\},
\end{align}
and, for continuous latent fields $v,w$, set pointwise
\begin{equation}\label{eq:q-prior-transform}
 q=L_{q_-,q_+}(v),\quad
 \gamma=\gamma_-(q)+[\gamma_+(q)-\gamma_-(q)]\ell(w),\quad \beta=\gamma+q.
\end{equation}
For $q\in(q_-,q_+)$, $\gamma_-(q)<\gamma_+(q)$, so $\beta$ and $\gamma$
lie strictly within their prescribed bounds. The map $(v,w)\mapsto
(\beta,\gamma)$ in \eqref{eq:q-prior-transform} is Lipschitz in the
supremum norm because $L_{q_-,q_+}$, $\ell$, $\gamma_-$, and $\gamma_+$
are Lipschitz and the width $\gamma_+-\gamma_-$ is bounded. This
map can induce dependence between the physical rates even when the
latent fields are independent. The bounds set the ranges of admissible
values. The means and covariances of the latent fields determine how
the prior distributes mass within those ranges. In the analysis we
fix these choices before conditioning on $y$. We refer to
\cref{app:latent-truncation} for a continuous Gaussian series construction
and a posterior error bound for its truncation.

\section{Low-prevalence asymptotics and weak identifiability}\label{sec:low-prevalence}
The decomposition in \eqref{eq:intro-I-rewrite} separates the linear
infected dynamics from the nonlinear saturation term. We use it to
expand the observation maps and to bound the effect of saturation on
the growth estimates.

\subsection{Uniform expansions}
Fix a nonnegative
$i_0\in L^\infty(\Omega)$ and, for $0<\varepsilon\le1$, replace the infected initial
condition by
\begin{equation}\label{eq:small-I0}
 I_0^\varepsilon=\varepsilon i_0.
\end{equation}
The susceptible initial condition remains fixed and satisfies $S_0\ge s_*>0$.
Write $(S_\theta^\varepsilon,I_\theta^\varepsilon)$ for the corresponding nonlinear
solution and set
\begin{equation}\label{eq:q-def}
 q_\theta=\beta-\gamma.
\end{equation}
Write $L_\theta=A_{d_I}+q_\theta$ on $L^2(\Omega)$, where $A_{d_I}$ is
defined in \cref{sec:heat-evolution} and $q_\theta$ acts by multiplication.
Let $i_\theta(t)=e^{tL_\theta}i_0$, so that
\begin{subequations}\label{eq:linear-I}
\begin{align}
 \partial_t i_\theta-d_I\Delta i_\theta
 &=q_\theta i_\theta &&\text{in }Q_T,\label{eq:linear-I-pde}\\
 \partial_\nu i_\theta&=0 &&\text{on }\partial\Omega\times(0,T),\\
 i_\theta(0)&=i_0 &&\text{in }\Omega.\label{eq:linear-I-ic}
\end{align}
\end{subequations}

\begin{theorem}\label{thm:low-prevalence}
There exists $C<\infty$, independent of $\theta\in\X_\ad$ and
$0<\varepsilon\le1$, such that
\begin{align}
 \norm{I_\theta^\varepsilon}_{L^\infty(Q_T)}
 &\le C\varepsilon,\label{eq:Ieps-small}\\
 \norm{I_\theta^\varepsilon-\varepsilon i_\theta}_{\Z_T^I}
 &\le C\varepsilon^2,\label{eq:Ieps-expansion}
\end{align}
where
\[
 \norm{v}_{\Z_T^I}
 =\norm{v}_{C([0,T];L^2)}+\norm{v}_{L^2(0,T;H^1)}.
\]
The same $O(\varepsilon^2)$ bound holds in $L^\infty(Q_T)$.
Consequently, if $\theta_1,\theta_2\in\X_\ad$ satisfy
\begin{equation}\label{eq:same-identifiable}
 q_{\theta_1}=q_{\theta_2},
 \qquad d_{I,1}=d_{I,2},
\end{equation}
then
\begin{equation}\label{eq:same-q-difference}
 \norm{I_{\theta_1}^\varepsilon-I_{\theta_2}^\varepsilon}_{\Z_T^I}
 \le C\varepsilon^2.
\end{equation}
\end{theorem}

\begin{proof}
Let $N_\theta^\varepsilon=S_\theta^\varepsilon+I_\theta^\varepsilon$.  Since
$\F(S,I)=I-I^2/(S+I)$, the infected equation is equivalent to
\begin{equation}\label{eq:I-exact-q}
 \partial_t I_\theta^\varepsilon-d_I\Delta I_\theta^\varepsilon
 =q_\theta I_\theta^\varepsilon
 -\beta\frac{(I_\theta^\varepsilon)^2}{N_\theta^\varepsilon}.
\end{equation}
The comparison argument in \cref{thm:forward-wp}, now using
$I_0^\varepsilon=\varepsilon i_0$, gives
\[
 \norm{I_\theta^\varepsilon(t)}_{L^\infty}
 \le\varepsilon e^{\overline\beta T}\norm{i_0}_{L^\infty},
\]
from which \eqref{eq:Ieps-small} follows.  Moreover,
$N_\theta^\varepsilon\ge S_\theta^\varepsilon\ge n_T$ by
\eqref{eq:N-lower}.

Set $R_\theta^\varepsilon=I_\theta^\varepsilon-\varepsilon i_\theta$.  Subtracting
$\varepsilon$ times \eqref{eq:linear-I} from \eqref{eq:I-exact-q} gives
\begin{equation}\label{eq:r-eqn}
 \partial_t R_\theta^\varepsilon-d_I\Delta R_\theta^\varepsilon
 =q_\theta R_\theta^\varepsilon
 -\beta\frac{(I_\theta^\varepsilon)^2}{N_\theta^\varepsilon},
 \qquad R_\theta^\varepsilon(0)=0.
\end{equation}
By \eqref{eq:Ieps-small} and the lower bound for $N_\theta^\varepsilon$,
\begin{equation}\label{eq:quadratic-force}
 \norm{\beta(I_\theta^\varepsilon)^2/N_\theta^\varepsilon}_{L^2(Q_T)}
 \le C\varepsilon^2.
\end{equation}
Testing \eqref{eq:r-eqn} with $R_\theta^\varepsilon$ and applying Young's
and Gronwall's inequalities gives \eqref{eq:Ieps-expansion}.
Recall $q_+=\overline\beta-\underline\gamma$ from \eqref{eq:q-feasible}.
For $v\in L^\infty(\Omega)$,
the function $u\mapsto e^{q_+u}\|v\|_{L^\infty}$ is a spatially constant
supersolution for both $e^{uL_\theta}v$ and $-e^{uL_\theta}v$. By the
parabolic comparison principle, as in \cref{thm:forward-wp},
\[
 \|e^{uL_\theta}v\|_{L^\infty}\le e^{q_+u}\|v\|_{L^\infty},
 \qquad u\ge0.
\]
The variation-of-constants formula from \cref{sec:heat-evolution}
applied to \eqref{eq:r-eqn} gives the $L^\infty(Q_T)$ remainder bound,
since the forcing is $O(\varepsilon^2)$
in $L^\infty(Q_T)$ by \eqref{eq:Ieps-small} and \eqref{eq:N-lower}. Under
\eqref{eq:same-identifiable}, the two leading terms
$\varepsilon i_{\theta_j}$ coincide.  Estimate \eqref{eq:same-q-difference} follows from the triangle inequality
and \eqref{eq:Ieps-expansion}.
\end{proof}

\begin{theorem}[Second-order expansion]\label{thm:second-order}
Let $S_\theta^{(0)}(t)=P_{d_S}(t)S_0$ be the Neumann heat evolution of
$S_0$, and let $j_\theta$ be the solution with $j_\theta(0)=0$ of
\begin{equation}\label{eq:second-order}
 \partial_t j_\theta-d_I\Delta j_\theta
 =q_\theta j_\theta-\frac{\beta i_\theta^2}{S_\theta^{(0)}},
 \qquad \partial_\nu j_\theta=0.
\end{equation}
There exists a constant $C<\infty$, independent of $\theta\in\X_\ad$ and
$0<\varepsilon\le1$, such that
\begin{equation}\label{eq:second-order-error}
 \|I_\theta^\varepsilon-\varepsilon i_\theta-\varepsilon^2j_\theta\|_{\Z_T^I}
 +\|I_\theta^\varepsilon-\varepsilon i_\theta-\varepsilon^2j_\theta\|_{L^\infty(Q_T)}
 \le C\varepsilon^3.
\end{equation}
In particular, the first correction to rescaled prevalence depends separately
on $\beta$, and may also depend on $d_S$ through $S_\theta^{(0)}$.
\end{theorem}
\begin{proof}
The absolute value of the reaction term in the susceptible equation is at
most $(\overline\beta+\overline\gamma)I_\theta^\varepsilon$. Applying
the variation-of-constants formula and $L^\infty$ contraction estimate
from \cref{sec:heat-evolution}, together with \eqref{eq:Ieps-small},
we obtain
$\|S_\theta^\varepsilon-S_\theta^{(0)}\|_{L^\infty(Q_T)}\le C\varepsilon$.
Let $N_\theta^\varepsilon=S_\theta^\varepsilon+I_\theta^\varepsilon$.
Since $S_\theta^{(0)}\ge s_*$ and $N_\theta^\varepsilon\ge n_T$, it follows that
\[
 \|N_\theta^\varepsilon-S_\theta^{(0)}\|_\infty\le C\varepsilon,
 \qquad
 \left\|\frac1{N_\theta^\varepsilon}
             -\frac1{S_\theta^{(0)}}\right\|_\infty\le C\varepsilon.
\]
The $L^\infty$ remainder in \cref{thm:low-prevalence} implies
$\|(I_\theta^\varepsilon)^2-\varepsilon^2 i_\theta^2\|_\infty
\le C\varepsilon^3$. Hence
\begin{equation}\label{eq:second-force-error}
 \left\|\frac{\beta(I_\theta^\varepsilon)^2}{N_\theta^\varepsilon}
       -\varepsilon^2\frac{\beta i_\theta^2}{S_\theta^{(0)}}\right\|_{L^\infty(Q_T)}
 \le C\varepsilon^3.
\end{equation}
Write $w=I_\theta^\varepsilon-\varepsilon i_\theta-\varepsilon^2j_\theta$.
Subtracting the equations for $\varepsilon i_\theta$ and
$\varepsilon^2j_\theta$ from \eqref{eq:I-exact-q}, we obtain
\[
 \partial_t w-d_I\Delta w=q_\theta w+g,\qquad
 \partial_\nu w=0,\qquad w(0)=0,
 \qquad
 g=\varepsilon^2\frac{\beta i_\theta^2}{S_\theta^{(0)}}
   -\frac{\beta(I_\theta^\varepsilon)^2}{N_\theta^\varepsilon}.
\]
By \eqref{eq:second-force-error}, $\|g\|_{L^\infty(Q_T)}\le C\varepsilon^3$.
Using the variation-of-constants formula and the estimate
$\|e^{uL_\theta}v\|_{L^\infty}\le e^{q_+u}\|v\|_{L^\infty}$ from the
proof of \cref{thm:low-prevalence}, we find that
\[
 \|w(t)\|_{L^\infty}
 \le \int_0^t e^{q_+(t-r)}\|g(r)\|_{L^\infty}\dd r
 \le C\varepsilon^3,\qquad 0\le t\le T.
\]
To estimate the $\Z_T^I$ norm, we use $w(t)\in H^1(\Omega)$ as a test
function in the weak equation for almost every $t$. The time derivative
contributes $\frac12\frac{\dd}{\dd t}\|w(t)\|_{L^2}^2$ in the
$H^1(\Omega)^*,H^1(\Omega)$ pairing. Young's inequality gives
\[
 \frac12\frac{\dd}{\dd t}\|w(t)\|_{L^2}^2+d_I\|\nabla w(t)\|_{L^2}^2
 \le \left(\|q_\theta\|_{L^\infty}+\frac12\right)\|w(t)\|_{L^2}^2
      +\frac12\|g(t)\|_{L^2}^2.
\]
Since $w(0)=0$ and $d_I\ge\underline d$, Gronwall's inequality and
integration give $\|w\|_{\Z_T^I}\le C\|g\|_{L^2(Q_T)}\le C\varepsilon^3$.
The constants depend only on $\Omega$, $T$, the coefficient bounds in
\eqref{eq:Xad-def}, the lower bound $s_*$, and the initial profiles
through $\|S_0\|_{L^\infty}$ and $\|i_0\|_{L^\infty}$.
\end{proof}

\subsection{Posterior limits and weak rate information}
We use a bounded linear prevalence observation operator
$\mathcal O_I:C([0,T];L^2(\Omega))\to\R^J$ as in \cref{sec:model}.
It also acts boundedly on $\Z_T^I$ because the $\Z_T^I$ norm controls
the $C([0,T];L^2)$ norm. The unscaled experiment becomes uninformative
as $\varepsilon\to0$ if the noise level is fixed. To retain a
nondegenerate signal-to-noise ratio, consider
\begin{equation}\label{eq:scaled-data}
 Y_\varepsilon
 =\mathcal O_I I_\theta^\varepsilon+\varepsilon\eta,
 \qquad \eta\sim N(0,\Gamma),
\end{equation}
or, equivalently,
\begin{equation}\label{eq:scaled-data-bar}
 \overline Y_\varepsilon
 =\varepsilon^{-1}Y_\varepsilon
 =\G_\varepsilon^{\prev}(\theta)+\eta,
 \qquad
 \G_\varepsilon^{\prev}(\theta)
 =\varepsilon^{-1}\mathcal O_I I_\theta^\varepsilon.
\end{equation}
The limiting map is
\begin{equation}\label{eq:G0-prev}
 \G_0^{\prev}(\theta)=\mathcal O_I i_\theta.
\end{equation}

\begin{theorem}[Prevalence observations]\label{thm:low-prevalence-posterior}
Let $\mu_\varepsilon^{\overline y}$ be the posterior associated with
\eqref{eq:scaled-data-bar}, and let $\mu_{\mathrm{lin}}^{\overline y}$ be defined by
\begin{equation}\label{eq:linear-posterior}
 \frac{\dd\mu_{\mathrm{lin}}^{\overline y}}{\dd\mu_0}(\theta)
 =\frac1{Z_{\mathrm{lin}}(\overline y)}
 \exp\left(
 -\frac12|\Gamma^{-1/2}
 (\overline y-\G_0^{\prev}(\theta))|^2
 \right).
\end{equation}
Then, for every $R<\infty$, there is $C_R<\infty$ such that
\begin{equation}\label{eq:low-prev-Hell}
 \sup_{|\overline y|\le R}
 d_\Hell(\mu_\varepsilon^{\overline y},
         \mu_{\mathrm{lin}}^{\overline y})
 \le C_R\varepsilon.
\end{equation}
The limiting likelihood in \eqref{eq:linear-posterior} depends on
$\theta=(\beta,\gamma,d_S,d_I)$ only through $(q_\theta,d_I)$.
\end{theorem}

\begin{proof}
By \cref{thm:low-prevalence}, uniformly over $\theta\in\X_\ad$,
\[
 |\G_\varepsilon^{\prev}(\theta)-\G_0^{\prev}(\theta)|
 \le \norm{\mathcal O_I}\,
 \varepsilon^{-1}
 \norm{I_\theta^\varepsilon-\varepsilon i_\theta}_{\Z_T^I}
 \le C\varepsilon.
\]
The linear parabolic energy estimate bounds $i_\theta$ uniformly in $\Z_T^I$.
Together with \eqref{eq:Ieps-expansion}, this also bounds
$I_\theta^\varepsilon/\varepsilon$ uniformly in $\Z_T^I$.
The boundedness of $\mathcal O_I$ then gives uniform bounds on both
observation maps. Apply \cref{thm:posterior-stability} with maps
$\G_\varepsilon^{\prev}$ and $\G_0^{\prev}$, and with
$y_1=y_2=\overline y$, to obtain \eqref{eq:low-prev-Hell}.
Equation \eqref{eq:linear-I} gives the asserted parameter dependence.
\end{proof}

\begin{remark}
\label{rem:incidence-information}
The conclusion of \cref{thm:low-prevalence-posterior} is specific to prevalence-type
observations.  Indeed,
\begin{equation}\label{eq:incidence-expansion}
 \beta\F(S_\theta^\varepsilon,I_\theta^\varepsilon)
 =\varepsilon\beta i_\theta+O(\varepsilon^2)
 \quad\text{in }C([0,T];L^2),
\end{equation}
uniformly over $\X_\ad$.  Thus rescaled incidence observations depend on $\beta$ at
first order, while the leading prevalence dynamics depend on $q=\beta-\gamma$.  Combining both
data types can therefore separate the transmission and recovery fields under conditions
where prevalence alone cannot.  Identifiability of the coefficient fields also depends on the number and
choice of observations.
\end{remark}

\begin{corollary}[Joint observations]
\label{cor:joint-limit}
Let $\mathcal O_I$ be the prevalence observation operator in
\eqref{eq:scaled-data}, and let
$\mathcal O_C:C([0,T];L^2(\Omega))\to\R^{J_C}$ be a bounded linear
incidence observation operator. Define
\begin{align*}
\G_\varepsilon^{\rm joint}(\theta)
 &=\bigl(\varepsilon^{-1}\mathcal O_I I_\theta^\varepsilon,
       \varepsilon^{-1}\mathcal O_C[\beta\F(S_\theta^\varepsilon,I_\theta^\varepsilon)]\bigr),\\
\G_{\rm lin}^{\rm joint}(\theta)
 &=\bigl(\mathcal O_I i_\theta,\mathcal O_C(\beta i_\theta)\bigr).
\end{align*}
With the same prior, datum $\overline y$, and fixed positive definite joint
noise covariance, the corresponding posteriors satisfy
\[
\sup_{|\overline y|\le R}
d_\Hell(\mu_\varepsilon^{\overline y,\rm joint},
\mu_{\rm lin}^{\overline y,\rm joint})\le C_R\varepsilon.
\]
\end{corollary}
\begin{proof}
We subtract $\varepsilon\beta i_\theta$ from the incidence identity
$\beta\F=\beta I-\beta I^2/(S+I)$. Theorem~\ref{thm:low-prevalence} and the
uniform denominator bound give an $O(\varepsilon^2)$ error in $C([0,T];L^2)$.
Thus the joint maps differ uniformly by $O(\varepsilon)$ and remain uniformly
bounded. Applying \cref{thm:posterior-stability} completes the proof.
\end{proof}

Two parameters belong to the same fibre of the map
$\theta\mapsto(q_\theta,d_I)$ if they satisfy \eqref{eq:same-identifiable}.

\begin{corollary}
\label{cor:KL-low-prevalence}\label{rem:weak-information}
Let $\overline P_{\varepsilon,\theta}
=N(\G_\varepsilon^{\prev}(\theta),\Gamma)$.
If $\theta_1,\theta_2$ satisfy \eqref{eq:same-identifiable}, then
\begin{equation}\label{eq:KL-order}
 D_{\rm KL}(\overline P_{\varepsilon,\theta_1}
       \Vert\overline P_{\varepsilon,\theta_2})
 =\frac{\varepsilon^2}{2}
   |\Gamma^{-1/2}\mathcal O_I(j_{\theta_1}-j_{\theta_2})|^2
   +O(\varepsilon^3).
\end{equation}
For fixed parameters with distinct limiting observation means, the divergence
instead tends to
$\frac12|\Gamma^{-1/2}(\G_0^{\prev}(\theta_1)-\G_0^{\prev}(\theta_2))|^2>0$.
\end{corollary}
\begin{proof}
For Gaussian laws with equal covariance, the KL divergence is half the
squared covariance-weighted difference of their means. We divide
\eqref{eq:second-order-error} by $\varepsilon$ and apply the bounded
observation operator. On a common fibre the leading means cancel, leaving
$\varepsilon\mathcal O_I(j_{\theta_1}-j_{\theta_2})+O(\varepsilon^2)$.
Equation \eqref{eq:KL-order} follows by squaring this expression. The
second claim follows from convergence of the observation means.
\end{proof}

We use \emph{weak rate information}, or weak identifiability here, to describe
this loss of sensitivity along the limiting fibres. The leading coefficient in
\eqref{eq:KL-order} is positive for an informative observation of the
second-order difference. For example, fix $q,d_I,d_S$ and increase both
rates by a nonzero continuous $h\ge0$, within the admissible bounds.
Then
\[
 j_{\theta_2}(u)-j_{\theta_1}(u)
 =-\int_0^u e^{(u-r)L_\theta}
       \frac{h\,i_\theta(r)^2}{S_\theta^{(0)}(r)}\dd r
 \qquad(u>0)
\]
where $L_\theta$, $i_\theta$ and $S_\theta^{(0)}$ are common to the two
parameters. For $v\ge0$, comparison with the constant potential $q_-$
gives $e^{tL_\theta}v\ge e^{q_-t}P_{d_I}(t)v$. By the positivity
improvement in \cref{sec:heat-evolution}, $i_\theta(r)>0$ for $r>0$
when $i_0\not\equiv0$. Thus the forcing in the integral is nonnegative
and nonzero for $0<r<u$, and the same comparison makes its propagated
value strictly positive. Consequently $j_{\theta_2}(u)-j_{\theta_1}(u)<0$
almost everywhere for $u>0$. A nonzero nonnegative spatial average at positive time
detects this difference.
General observation operators can annihilate it. At fixed unscaled noise,
the corresponding fibre divergence is $O(\varepsilon^4)$, while arbitrary
fixed alternatives already have divergence $O(\varepsilon^2)$. The rescaled
experiment therefore distinguishes these two orders of sensitivity.

\subsection{Prior dependence along likelihood fibres}\label{sec:prior-fibres}
The next proposition specializes the conditional-prior principle for Bayesian
nonidentified models \cite{Poirier1998,Gustafson2015} to the limiting SIS likelihood.
\begin{proposition}
\label{prop:conditional-prior}
Let $H(\theta)=(\beta-\gamma,d_I)$, $\pi_0=H_\sharp\mu_0$, and
$\pi_{\mathrm{lin}}^{\overline y}=H_\sharp\mu_{\mathrm{lin}}^{\overline y}$.
There exists a regular conditional prior $\mu_0(\dd\theta\mid h)$, and the
limiting prevalence posterior satisfies
\begin{equation}\label{eq:conditional-prior}
 \mu_{\mathrm{lin}}^{\overline y}(\dd\theta)
 =\int \mu_0(\dd\theta\mid h)\,
       \pi_{\mathrm{lin}}^{\overline y}(\dd h),\qquad
 \frac{\dd\pi_{\mathrm{lin}}^{\overline y}}{\dd\pi_0}(h)
 =\frac{e^{-\Phi_{\mathrm{lin}}(h;\overline y)}}{Z_{\mathrm{lin}}(\overline y)},
\end{equation}
where $\Phi_{\mathrm{lin}}$ is the Gaussian potential for \eqref{eq:G0-prev} and
$Z_{\mathrm{lin}}(\overline y)=\int e^{-\Phi_{\mathrm{lin}}(h;\overline y)}\pi_0(\dd h)$.
Consequently, the conditional posterior given $H=h$ equals the conditional prior
for $\pi_0$-almost every $h$.
\end{proposition}

\begin{proof}
The admissible set is a closed subset of a separable Banach space, and $H$ is
continuous with values in the separable Banach space
$\mathcal H=C(\overline\Omega)\times\R$. A regular conditional prior therefore
exists. For Borel sets $A\subset\X_\ad$ and $B\subset\mathcal H$, it satisfies
\[
 \mu_0(A\cap H^{-1}(B))=\int_B\mu_0(A\mid h)\,\pi_0(\dd h).
\]
Write $\ell_{\overline y}(h)=e^{-\Phi_{\mathrm{lin}}(h;\overline y)}$, so that
the likelihood in \eqref{eq:linear-posterior} is $\ell_{\overline y}(H(\theta))$.
Applying Bayes' formula and the conditional prior identity gives
\[
 \mu_{\mathrm{lin}}^{\overline y}(A\cap H^{-1}(B))
 =\frac{\displaystyle\int_B\mu_0(A\mid h)\ell_{\overline y}(h)\,\pi_0(\dd h)}
        {\displaystyle\int_{\mathcal H}\ell_{\overline y}(h)\,\pi_0(\dd h)}.
\]
Taking $A=\X_\ad$ gives the marginal density in \eqref{eq:conditional-prior}.
The denominator equals $Z_{\mathrm{lin}}(\overline y)$ by the definition of
$\pi_0=H_\sharp\mu_0$ and makes this density integrate to one.
The displayed identity then shows that
$\mu_0(\dd\theta\mid h)$ is also a version of the conditional posterior.
Taking $B=\mathcal H$ gives the mixture in \eqref{eq:conditional-prior}.
Since $\ell_{\overline y}>0$, the marginal posterior and $\pi_0$ are equivalent,
so the conditional identity holds for $\pi_0$-almost every $h$.
\end{proof}

Thus marginal recovery or transmission uncertainty can change through prior
dependence on $H$, even though the limiting likelihood supplies no additional
information within a fibre.  For a numerical study, this suggests reporting both $q$ and the separate rates, and comparing priors
with the same $H$-marginal but different conditional rate distributions.  For the
limiting experiment these priors give the same posterior law of $H$, and hence of
the persistence event, while the separate rate laws may differ.

\subsection{From prevalence growth to the persistence threshold}
\label{sec:growth-threshold}

The linearised prevalence growth rate depends on the same pair $(q,d_I)$
as the threshold eigenvalue. We relate the finite-time growth rate to this
eigenvalue and then account for nonlinear saturation. The operator $L_\theta$
is self-adjoint on $L^2(\Omega)$. Denote its eigenvalues,
including multiplicities, by
\[
 \lambda_1(\theta)=\lambda_*(\theta)>
 \lambda_2(\theta)\ge\lambda_3(\theta)\ge\cdots,
\]
and choose an orthonormal eigenbasis $\{\phi_k\}_{k\ge1}$ with
$\phi_1>0$.  For the nonnegative, nonzero initial profile $i_0$, set
\begin{equation}\label{eq:spectral-transient-constants}
 c_k=\ip{i_0}{\phi_k}_{L^2},\qquad
 \mathfrak g_\theta=\lambda_1-\lambda_2,
 \qquad
 \rho_\theta^2=\frac{\sum_{k\ge2}c_k^2}{c_1^2}.
\end{equation}
The strong positivity of $\phi_1$ gives $c_1>0$, while simplicity of the principal
eigenvalue gives $\mathfrak g_\theta>0$.

For $0\le s<t\le T$, define the finite-time prevalence growth rates
\begin{align}
 \Lambda_\theta(s,t)
 &=\frac1{t-s}\log\frac{\norm{i_\theta(t)}_{L^2}}
                              {\norm{i_\theta(s)}_{L^2}},
 \label{eq:linear-growth-rate}\\
 \Lambda_\theta^\varepsilon(s,t)
 &=\frac1{t-s}\log\frac{\norm{I_\theta^\varepsilon(t)}_{L^2}}
                              {\norm{I_\theta^\varepsilon(s)}_{L^2}}.
 \label{eq:nonlinear-growth-rate}
\end{align}
The factor $\varepsilon$ cancels in the second ratio.

\begin{theorem}[Nonlinear growth bounds]\label{thm:growth-threshold}
For $\theta\in\X_\ad$, $i_0\ge0$, $i_0\not\equiv0$, and $0\le s<t\le T$,
the linear growth rate satisfies
\begin{equation}\label{eq:linear-growth-bracket}
 \Lambda_\theta(s,t)\le\lambda_*(\theta)
 \le\Lambda_\theta(s,t)+\mathcal E_\theta(s,t),\qquad
 \mathcal E_\theta(s,t)=
 \frac{\log(1+\rho_\theta^2e^{-2\mathfrak g_\theta s})}{2(t-s)}.
\end{equation}
Let $U\in L^1(0,t)$ be nonnegative and satisfy
$\|I_\theta^\varepsilon(r)\|_\infty\le U(r)$ for almost every $r\in(0,t)$.
Set
\begin{equation}\label{eq:growth-nonlinear-error}
 A_U(u)=\frac{\overline\beta}{s_*}
          \int_0^u e^{\overline\beta r}U(r)\dd r,
 \qquad B_U(s,t)=\frac{A_U(t)}{t-s}.
\end{equation}
The nonlinear growth satisfies
\begin{equation}\label{eq:nonlinear-growth-bracket}
 \Lambda_\theta^\varepsilon(s,t)\le\lambda_*(\theta)
 \le\Lambda_\theta^\varepsilon(s,t)+\mathcal E_\theta(s,t)+B_U(s,t).
\end{equation}
The lower bound holds at every infection size, independently of $U$.

Suppose measured norms $\widehat N_s,\widehat N_t>0$ satisfy
\begin{equation}\label{eq:relative-norm-error}
 \left|\frac{\widehat N_u}{\|I_\theta^\varepsilon(u)\|_{L^2}}-1\right|
 \le\delta_{\rm obs}<1,\qquad u\in\{s,t\}.
\end{equation}
Define
\[
 \widehat\Lambda=\frac1{t-s}\log\frac{\widehat N_t}{\widehat N_s},\qquad
 d_{\rm obs}=\frac1{t-s}\log\frac{1+\delta_{\rm obs}}{1-\delta_{\rm obs}}.
\]
Then the measured growth rate gives the bounds
\begin{equation}\label{eq:measured-growth-bracket}
 \widehat\Lambda-d_{\rm obs}\le\lambda_*(\theta)
 \le\widehat\Lambda+\mathcal E_\theta+B_U+d_{\rm obs}.
\end{equation}
In particular, we have
\begin{align}
 \widehat\Lambda>d_{\rm obs}
 &\ \Longrightarrow\ \Rzero(\theta)>1,\label{eq:persistence-certificate}\\
 \widehat\Lambda+\mathcal E_\theta+B_U+d_{\rm obs}<0
 &\ \Longrightarrow\ \Rzero(\theta)<1.\label{eq:extinction-certificate}
\end{align}
\end{theorem}
\begin{proof}
The spectral expansion of $L_\theta$ gives
$i_\theta(u)=\sum_{k\ge1}c_ke^{\lambda_ku}\phi_k$. By orthogonality,
\begin{equation}\label{eq:growth-spectral-expansion}
 \|i_\theta(u)\|_{L^2}^2=c_1^2e^{2\lambda_1u}(1+r_\theta(u)),\qquad
 r_\theta(u)=\sum_{k\ge2}\frac{c_k^2}{c_1^2}
                  e^{-2(\lambda_1-\lambda_k)u}.
\end{equation}
Since $0\le r_\theta(t)\le r_\theta(s)\le
\rho_\theta^2e^{-2\mathfrak g_\theta s}$,
\[
 \lambda_1-\Lambda_\theta(s,t)
 =\frac{\log(1+r_\theta(s))-\log(1+r_\theta(t))}{2(t-s)}
 \in[0,\mathcal E_\theta(s,t)].
\]
The linear bound \eqref{eq:linear-growth-bracket} follows.

For the nonlinear lower bound, write $I=I_\theta^\varepsilon$ and
$M_I(u)=\int_\Omega I(u)\dd x$. Taking the test function one in the
weak infected equation gives $M_I'(u)\ge-\overline\gamma M_I(u)$,
since the incidence is nonnegative. Hence
$M_I(u)\ge e^{-\overline\gamma u}M_I(0)>0$ for all $u$, and
$\|I(u)\|_{L^2}>0$. Testing \eqref{eq:I-exact-q} with $I$ and using
\eqref{eq:lambda-def} gives
\begin{equation}\label{eq:exact-nonlinear-energy}
 \frac12\frac{\dd}{\dd u}\|I(u)\|_{L^2}^2
 =-d_I\|\nabla I\|_{L^2}^2+\int_\Omega q I^2\dd x
       -\int_\Omega\frac{\beta I^3}{S+I}\dd x
 \le\lambda_*\|I(u)\|_{L^2}^2.
\end{equation}
The energy chain rule holds for the weak solution. Integrating the
inequality gives $\|I(t)\|_{L^2}\le e^{\lambda_*(t-s)}\|I(s)\|_{L^2}$.
Taking the logarithmic ratio yields $\Lambda_\theta^\varepsilon\le\lambda_*$.

For the upper bound, we write the equation as
$\partial_u I=L_\theta I-\delta_\theta I$, where
$\delta_\theta=\beta I/(S+I)$. The susceptible lower bound implies
\[
 0\le\delta_\theta(x,u)
 \le \overline\beta e^{\overline\beta u}U(u)/s_*=A_U'(u).
\]
Comparison with the linear solution therefore gives
\begin{equation}\label{eq:multiplicative-comparison}
 e^{-A_U(u)}\varepsilon i_\theta(u)
 \le I_\theta^\varepsilon(u)\le\varepsilon i_\theta(u),\qquad 0\le u\le t.
\end{equation}
Indeed, the left member is a subsolution of the equation with coefficient
$\delta_\theta$, and the right member is a supersolution, with identical
initial data. Applying this comparison at $s,t$ gives
$\Lambda_\theta(s,t)\le\Lambda_\theta^\varepsilon(s,t)+A_U(t)/(t-s)$.
Combining this inequality with \eqref{eq:linear-growth-bracket} establishes
\eqref{eq:nonlinear-growth-bracket}.
Finally, \eqref{eq:relative-norm-error} implies
$|\widehat\Lambda-\Lambda_\theta^\varepsilon|\le d_{\rm obs}$.
We therefore obtain \eqref{eq:measured-growth-bracket}. The sufficient
conditions for persistence and extinction follow from the threshold
equivalence in \cref{thm:threshold-equivalence}.
\end{proof}

An available uniform envelope is
$U(u)=\|I_0^\varepsilon\|_\infty e^{q_+u}$,
$q_+=\overline\beta-\underline\gamma$, because
$\partial_t I-d_I\Delta I\le q_+I$. It gives the explicit width
\begin{equation}\label{eq:explicit-nonlinear-width}
 B_U(s,t)=\frac{\overline\beta\|I_0^\varepsilon\|_\infty}{s_*(t-s)}
                 \int_0^t e^{(\overline\beta+q_+)u}\dd u.
\end{equation}
For $I_0^\varepsilon=\varepsilon i_0$ this is $O(\varepsilon)$ at fixed
$s,t$. It can be evaluated from an upper bound on the absolute initial
infection, without identifying a separate normalization $\varepsilon$.
A verified trajectory envelope can replace this conservative bound.
Use of a trajectory envelope requires an independent bound on the spatial
maximum of the infection over the observation period.

\begin{proposition}\label{prop:uniform-spectral}
Let $V=|\Omega|$, $q_-=\underline\beta-\overline\gamma$, and let
$\nu_1>0$ be the first positive Neumann eigenvalue of $-\Delta$.
There is $\mathfrak h\ge1$, depending only on $\Omega$, the coefficient bounds, and
a fixed heat-kernel comparison time, such that uniformly over $\X_\ad$ we have
\begin{align}
 \mathfrak g_\theta&\ge g_0:=\underline d\,\nu_1/\mathfrak h^2,\notag\\
 \rho_\theta^2&\le \rho_{\rm max}^2:=
 V\mathfrak h^2\|i_0\|_{L^2}^2/\|i_0\|_{L^1}^2-1.\label{eq:uniform-spectral}
\end{align}
Consequently, the transient term satisfies
\begin{equation}\label{eq:uniform-transient}
 \mathcal E_\theta(s,t)\le\overline{\mathcal E}(s,t)
 :=\frac{\log(1+\rho_{\rm max}^2e^{-2g_0s})}{2(t-s)}.
\end{equation}
\end{proposition}
\begin{proof}
The heat-semigroup comparison in \cref{app:spectral} gives
$\max\phi_1/\min\phi_1\le \mathfrak h$ for the $L^2$-normalized principal
eigenfunction. Hence $c_1\ge\|i_0\|_{L^1}/(\mathfrak h\sqrt V)$, and the
bound on $\rho_\theta$ follows from Parseval's identity. The ground-state transform
expresses the spectral gap as the weighted Poincar\'e quotient with weight
$\phi_1^2$. Comparison with the unweighted Neumann Poincar\'e inequality
gives $\mathfrak g_\theta\ge d_I\nu_1/\mathfrak h^2$.
\Cref{app:spectral} supplies the kernel constants and the full argument.
\end{proof}

Using \eqref{eq:explicit-nonlinear-width} and \eqref{eq:uniform-transient}
makes the error widths independent of the unknown parameter.
Numerical evaluation of these sufficient conditions requires verified bounds
on $\mathfrak h$ and $\nu_1$. An interval-domain construction is given in \cref{app:spectral}.
For positive weighted counts $\langle I,w\rangle$, including total infected
counts, \cref{prop:weighted-growth} gives a two-sided variant. Its signed
spectral remainder lacks the exact lower bound enjoyed by $L^2$ norms.

\begin{corollary}\label{cor:posterior-growth-bracket}
Let $\mu^{y,\varepsilon}$ be a nonlinear posterior, and suppose envelopes
$U_\theta$ are available for $\mu^{y,\varepsilon}$-almost every parameter.
Then the posterior persistence probability satisfies
\begin{align}
 \mu^{y,\varepsilon}(\Lambda_\theta^\varepsilon>0)
 &\le\mu^{y,\varepsilon}(\lambda_*>0)\notag\\
 &\le\mu^{y,\varepsilon}
       (\Lambda_\theta^\varepsilon+\mathcal E_\theta+B_{U_\theta}>0).
 \label{eq:posterior-growth-bracket}
\end{align}
For the linear posterior set $B_{U_\theta}=0$ and replace
$\Lambda_\theta^\varepsilon$ by $\Lambda_\theta$.
\end{corollary}
\begin{proof}
\Cref{thm:growth-threshold} gives the corresponding pointwise event
inclusions. The claim follows by taking posterior probabilities.
\end{proof}

These brackets interpret threshold uncertainty through growth evidence.
When a posterior parameter sample is already available, directly solving
its threshold eigenproblem remains the natural way to evaluate persistence.
The data-level certificates are useful when independent bounds on their
error terms are available.

\section{Posterior uncertainty in the epidemic threshold}\label{sec:threshold}

The following estimates control posterior propagation and approximation of
the threshold quantities defined in \cref{sec:threshold-definitions}.

\begin{theorem}
\label{thm:threshold-stability}
There is $C<\infty$, depending only on $\X_\ad$ and $\Omega$, such that
\begin{align}
 |\lambda_*(\theta_1)-\lambda_*(\theta_2)|
 &\le C\left(
 \norm{(\beta_1-\gamma_1)-(\beta_2-\gamma_2)}_{L^\infty}
 +|d_{I,1}-d_{I,2}|
 \right),\label{eq:lambda-Lipschitz}\\
 |\Rzero(\theta_1)-\Rzero(\theta_2)|
 &\le C\left(
 \norm{\beta_1-\beta_2}_{L^\infty}
 +\norm{\gamma_1-\gamma_2}_{L^\infty}
 +|d_{I,1}-d_{I,2}|
 \right).
 \label{eq:R0-Lipschitz}
\end{align}
\end{theorem}

\begin{proof}
Write $q_j=\beta_j-\gamma_j$ and let $\phi_j$ be an $L^2$-normalized maximizer for
$\lambda_*(q_j,d_{I,j})$.  The variational characterization gives
\begin{align*}
 \lambda_*(q_1,d_{I,1})-\lambda_*(q_2,d_{I,2})
 &\le\int_\Omega(q_1-q_2)\phi_1^2\dd x
 -(d_{I,1}-d_{I,2})\int_\Omega|\nabla\phi_1|^2\dd x.
\end{align*}
Since $|q_j|\le\overline\beta+\overline\gamma=:M_q$ and
$\lambda_*(q_j,d_{I,j})\ge-M_q$, the eigenvalue identity implies
\[
 \int_\Omega|\nabla\phi_j|^2\dd x
 \le\frac{2M_q}{\underline d}.
\]
Interchanging the two parameters, we obtain \eqref{eq:lambda-Lipschitz}.

For $\Rzero$, define $A_\theta(\varphi)=\int_\Omega\beta\varphi^2\dd x$
and $B_\theta(\varphi)=d_I\int_\Omega|\nabla\varphi|^2\dd x+\int_\Omega\gamma\varphi^2\dd x$.
Uniform coercivity of $B_\theta$ and boundedness of $A_\theta$ show that
$A_\theta/B_\theta$ changes by at most the right-hand side of
\eqref{eq:R0-Lipschitz}, uniformly over nonzero $\varphi$.  Taking suprema over $\varphi$ completes the proof.
\end{proof}

\begin{remark}
The threshold sign in \eqref{eq:threshold-equivalence} depends only on $(q,d_I)$, the
same pair that appears in the limiting prevalence likelihood of
\cref{thm:low-prevalence-posterior}.  The numerical value of $\Rzero$, in contrast,
depends separately on $\beta$ and $\gamma$ through \eqref{eq:R0-def}.  Hence early
prevalence data may identify the persistence event more reliably than the numerical
value of $\Rzero$ or the individual rate fields.
\end{remark}

\subsection{Posterior laws and persistence probabilities}

For a Borel map $Q:\X_\ad\to\R$, write $Q_\sharp\mu$ for the pushforward of $\mu$.

\begin{theorem}\label{thm:threshold-posterior}
Let $Q$ be either $\Rzero$ or $\lambda_*$.  Then $Q$ is bounded and Borel measurable on
$\X_\ad$. For every $R<\infty$ and $|y_1|,|y_2|\le R$, the pushforward laws satisfy
\begin{equation}\label{eq:pushforward-Hell}
 d_\Hell(Q_\sharp\mu^{y_1},Q_\sharp\mu^{y_2})
 \le d_\Hell(\mu^{y_1},\mu^{y_2})
 \le C_R|y_1-y_2|.
\end{equation}
For every $c\in\R$, we also have
\begin{align}
 \left|\E^{\mu^{y_1}}Q-\E^{\mu^{y_2}}Q\right|
 &\le2\sqrt2\norm{Q}_{L^\infty(\X_\ad)}
 d_\Hell(\mu^{y_1},\mu^{y_2}),\label{eq:Qmean-stability}\\
 \left|\mu^{y_1}(Q>c)-\mu^{y_2}(Q>c)\right|
 &\le\sqrt2d_\Hell(\mu^{y_1},\mu^{y_2}).
 \label{eq:Qprob-stability}
\end{align}
In particular, the posterior persistence probability
\begin{equation}\label{eq:persistence-prob}
 \Pp(\Rzero>1\mid y)=\mu^y(\Rzero>1)
 =\mu^y(\lambda_*>0)
\end{equation}
is locally Lipschitz stable with respect to the data.
\end{theorem}

\begin{proof}
Boundedness and measurability follow from \cref{thm:threshold-stability} and the bounded
admissible set.  Let $f$ and $g$ be the two posterior densities with respect to the prior
$\mu_0$. The corresponding densities of the laws of $Q$ with respect to
its prior law, evaluated at $Q$, are $\E_{\mu_0}[f\mid Q]$ and
$\E_{\mu_0}[g\mid Q]$. By conditional Cauchy--Schwarz,
\[
\int \sqrt{\E_{\mu_0}[f\mid Q]\,\E_{\mu_0}[g\mid Q]}\,\mathrm d\mu_0
\ge
\int \E_{\mu_0}[\sqrt{fg}\mid Q]\,\mathrm d\mu_0
=
\int \sqrt{fg}\,\mathrm d\mu_0.
\]
Subtracting each side from one gives the comparison of the squared
Hellinger distances and hence the first inequality in
\eqref{eq:pushforward-Hell}. The second inequality follows from
\eqref{eq:data-Hell-Lipschitz}.

The difference of the posterior means is bounded by
$2\|Q\|_\infty d_\TV$, and the difference of the probabilities of a
threshold event is bounded by $d_\TV$. Using
$d_\TV\le\sqrt{2}d_\Hell$ gives \eqref{eq:Qmean-stability} and
\eqref{eq:Qprob-stability}. Finally, the two events in
\eqref{eq:persistence-prob} coincide by \eqref{eq:threshold-equivalence}.
\end{proof}

We next account for simultaneous approximation of the posterior and of the threshold
functional.  This matters because $\Rzero$ is itself computed through an elliptic
eigenvalue problem.

\begin{theorem}\label{thm:threshold-approx}
Let $\mu_{\mathfrak a}^y$ satisfy \eqref{eq:posterior-approx-rate}, and let
$Q_{\mathfrak a}:\X_\ad\to\R$ be measurable with
\begin{equation}\label{eq:Qa-error}
 \sup_{\theta\in\X_\ad}|Q(\theta)-Q_{\mathfrak a}(\theta)|\le\delta_{\mathfrak a},
 \qquad \delta_{\mathfrak a}\to0,
\end{equation}
where $Q=\Rzero$ or $Q=\lambda_*$. The posterior mean satisfies
\begin{equation}\label{eq:threshold-mean-approx}
 \left|\int Q\dd\mu^y-\int Q_{\mathfrak a}\dd\mu_{\mathfrak a}^y\right|
 \le2\sqrt2\norm{Q}_{L^\infty}
 d_\Hell(\mu^y,\mu_{\mathfrak a}^y)+\delta_{\mathfrak a}.
\end{equation}
For any threshold $c\in\R$, the probability error satisfies
\begin{align}
 \left|\mu^y(Q>c)-\mu_{\mathfrak a}^y(Q_{\mathfrak a}>c)\right|
 &\le\sqrt2d_\Hell(\mu^y,\mu_{\mathfrak a}^y)
 +\mu_{\mathfrak a}^y(|Q-c|\le\delta_{\mathfrak a}).
 \label{eq:threshold-prob-approx}
\end{align}
Consequently, if $d_\Hell(\mu^y,\mu_{\mathfrak a}^y)\to0$ and
$\mu^y(Q=c)=0$, then
\begin{equation}\label{eq:threshold-prob-convergence}
 \mu_{\mathfrak a}^y(Q_{\mathfrak a}>c)\longrightarrow\mu^y(Q>c).
\end{equation}
\end{theorem}

\begin{proof}
Estimate \eqref{eq:threshold-mean-approx} follows by adding and subtracting
$\int Q\dd\mu_{\mathfrak a}^y$ and using the total-variation bound from
the proof of \cref{thm:posterior-approximation}.
For the probability, we write
\begin{align*}
 |\mu^y(Q>c)-\mu_{\mathfrak a}^y(Q_{\mathfrak a}>c)|
 \le|\mu^y(Q>c)-\mu_{\mathfrak a}^y(Q>c)| +\mu_{\mathfrak a}^y\bigl(\{Q>c\}\mathbin\triangle\{Q_{\mathfrak a}>c\}\bigr).
\end{align*}
The first term is bounded by total variation.  Under \eqref{eq:Qa-error}, the symmetric
difference is contained in $\{|Q-c|\le\delta_{\mathfrak a}\}$, so \eqref{eq:threshold-prob-approx} follows.  If $\mu^y(Q=c)=0$, then
$\mu^y(|Q-c|\le\delta)\to0$ as $\delta\downarrow0$.  Total-variation convergence of
$\mu_{\mathfrak a}^y$ to $\mu^y$ transfers the same property to the shrinking-band term, from which
\eqref{eq:threshold-prob-convergence} follows.
\end{proof}

Conforming eigenvalue approximation provides a concrete
choice of $Q_{\mathfrak a}$: quadratic Rayleigh-quotient error gives
$0\le\Rzero-\Rzero{}_h\le Ch^{2r}$ when the normalized principal eigenfunction
has uniform $H^1$ approximation error $O(h^r)$ \cite{Boffi2010}.
The precise approximation hypothesis and its scope are stated in \cref{cor:eigen-FE}.

\begin{remark}\label{rem:prior-no-atom}
The construction \eqref{eq:q-prior-transform} also provides a sufficient condition
for the no-threshold-atom assumption used in \cref{thm:threshold-approx}.
Suppose $v=c+\widetilde v$, where $c$ is a scalar random variable with an atomless
law, independent of $(\widetilde v,d_I)$. For instance, $c$ can be the
constant Gaussian mode in \eqref{eq:cosine-prior} with independent $d_I$.
Conditional on $(\widetilde v,d_I)$, increasing $c$ strictly increases $q$ at every
point.  For $c_2>c_1$, compactness of $\overline\Omega$ gives
$\min_{\overline\Omega}(q_{c_2}-q_{c_1})>0$, and \eqref{eq:lambda-def} implies
\[
 \lambda_*(q_{c_2},d_I)-\lambda_*(q_{c_1},d_I)
 \ge\min_{\overline\Omega}(q_{c_2}-q_{c_1})>0.
\]
There is at most one value of $c$ with $\lambda_*=0$.  Conditioning and using
atomlessness therefore give $\mu_0(\lambda_*=0)=0$.  Since the Gaussian-data
posterior is equivalent to its prior, $\mu^y(\lambda_*=0)=0$, and the threshold
equivalence gives $\mu^y(\Rzero=1)=0$ as well.
\end{remark}

\section{Discussion and numerical illustrations}\label{sec:numerics}\label{sec:conclusion}

The leading prevalence dynamics and the persistence event depend on the
same pair $(q,d_I)$. As \cref{prop:conditional-prior} shows, the limiting
prevalence likelihood updates the law of $(q,d_I)$ while preserving
the conditional prior on the separate rates given this pair. A posterior
can therefore assign high probability to one side of the threshold even when
transmission and recovery remain uncertain. Incidence can reduce that
uncertainty because its leading term depends explicitly on $\beta$.

The growth estimates also relate persistence to changes in prevalence
between two observation times. A measured growth rate of the infected
$L^2$ norm above the error allowance in \eqref{eq:measured-growth-bracket}
is sufficient for persistence. A decline can also occur with
$\lambda_*>0$ because of spatial transients or nonlinear saturation.
The upper bound in \eqref{eq:measured-growth-bracket} specifies how large
a decline is sufficient to infer extinction. Using this bound with
sparse measurements requires control of the error in recovering the
profile norm. For regional totals, \cref{prop:weighted-growth} gives
an estimate based on weighted integrals of the infected density.

To illustrate the posterior behaviour, we take $\Omega=(0,1)$,
$d_I=0.02$, $d_S=0.04$, $S_0=1$, and
$i_0=1+0.15\cos(\pi x)$, and parameterize the rates by
\begin{equation}\label{eq:numerical-parameterization}
 q=q_0+q_1\cos(\pi x),\qquad\gamma=a,\qquad\beta=a+q.
\end{equation}
We generate data at $(q_0,q_1,a)=(-0.022,0.120,0.340)$, for which the
inversion grid gives $\lambda_*^\dagger\approx\NumTrueLambda$. We observe
eleven cell averages at each of $t=0.5$ and $t=2$. We use independent Gaussian errors with
standard deviations $0.040$ for rescaled prevalence and $0.025$ for
rescaled incidence. We compute the posteriors by quadrature over the
three parameters $(q_0,q_1,a)$. For details of priors, observations and solvers, see \cref{app:implementation}.

\begin{figure}[htbp]
 \centering
 \includegraphics[width=\textwidth]{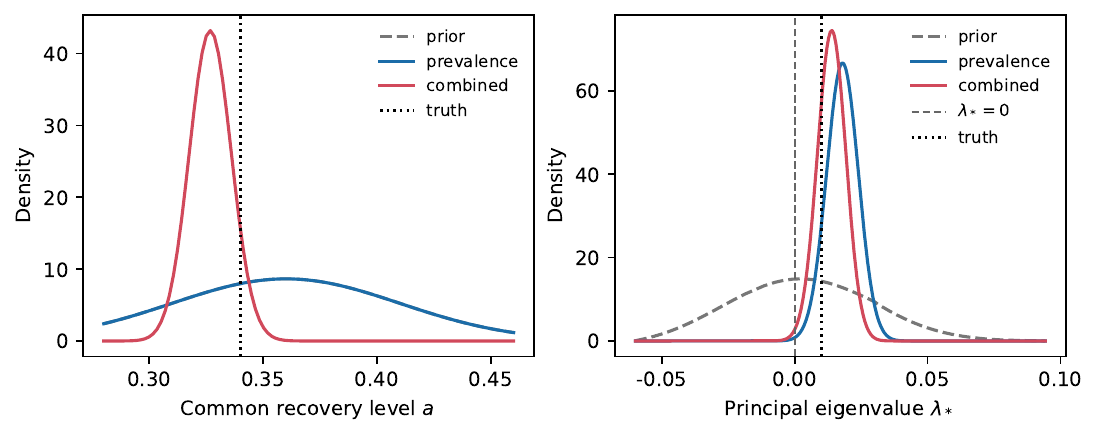}
 \caption{Controlled limiting-model experiment. The left panel shows
 the marginal densities of the recovery level $a$. The prevalence posterior
 and prior curves coincide. The right panel shows the corresponding
 densities of $\lambda_*$. Vertical lines mark the true values and, on the
 right, zero.}
 \label{fig:posterior-results}
\end{figure}

In the first experiment, we use the same limiting model and spatial grid
for data generation and inversion. Since $a$ is independent of $(q_0,q_1)$
under the prior, \cref{prop:conditional-prior} gives the same marginal
law of $a$ before and after conditioning on prevalence. The curves
coincide in the left panel of \cref{fig:posterior-results}, with a 95\% equal-tailed
interval of $[\NumPrevALow,\NumPrevAHigh]$. Adding incidence makes the
likelihood depend on $a$ and narrows the interval to
$[\NumBothALow,\NumBothAHigh]$. With dependent priors, a prevalence update
of $(q_0,q_1)$ can also change the marginal of $a$. The right panel shows
that prevalence alone concentrates the posterior on positive values of
$\lambda_*$. The persistence probability rises from
$\NumPriorThresholdProb$ to $\NumPrevThresholdProb$, and is
$\NumBothThresholdProb$ with both kinds of data.

For the discretized linear model, \cref{cor:posterior-growth-bracket}
gives the bracket $[\NumPrevGrowthLower,\NumPrevGrowthUpper]$ for the
posterior persistence probability. The lower endpoint is the probability
that the infected $L^2$ norm grows. With the prior support bounds $-0.24\le q\le0.21$ and $\tau_0=6$, the
interval estimate in \cref{app:spectral} gives the uniform transient
allowance $\overline{\mathcal E}(0.5,2)\approx2.79$.
The smallest linear growth rate on the grid is about $-0.061$, so with this
allowance the extinction condition is inconclusive even for exact norms.

\begin{figure}[!htbp]
 \centering
 \includegraphics[width=0.95\textwidth]{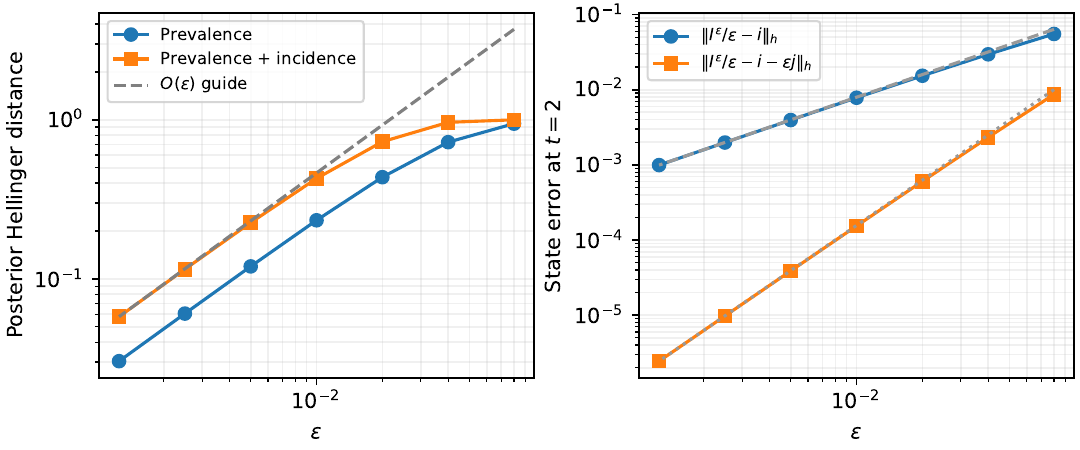}
 \caption{The left panel shows Hellinger distances between the nonlinear
 and limiting posteriors, conditioned on the same nonlinear data at each
 $\varepsilon$. The dashed line indicates order one. The right panel shows
 discrete $L^2$ errors of the rescaled state at $t=2$ on 41 cells, before and
 after the second-order correction. The dashed and dotted lines indicate
 orders one and two, respectively.}
 \label{fig:posterior-limit}
\end{figure}

To assess the nonlinear correction, we next generate data from the full
system on 123 cells and evaluate both likelihoods on 41 cells. At each
$\varepsilon$ we use the same rescaled data, prior, and covariance for
the two posteriors. We halve $\varepsilon$ from $0.080$ to $0.00125$.
The fitted Hellinger slopes
in \cref{fig:posterior-limit} over the three smallest scales are
\NumNonlinearPrevSlope\ for prevalence and \NumNonlinearJointSlope\
for joint data, in agreement with the first-order posterior bounds.
At $\varepsilon=0.01$ the distances are still $0.23$ and $0.43$, respectively,
and at $\varepsilon=0.005$ they are $0.12$ and $0.23$.
The fitted slopes for the rescaled state errors are \NumFirstOrderSlope\
before the correction and \NumSecondOrderSlope\ after it. We report
refinement checks in \cref{app:nonlinear-posterior}.
\FloatBarrier

\medskip
\begingroup
\clubpenalty=10000
\noindent\textbf{Concluding remarks.}\enspace
The posterior limits in \cref{thm:low-prevalence-posterior,cor:joint-limit}
describe a shrinking initial infection at a fixed signal-to-noise scale.
It would be of interest to study contraction of the posterior for
$\lambda_*$ under increasingly informative data, for example from repeated
measurements or from noise that decreases relative to the infection size.
A first question is when prevalence and incidence observations distinguish
parameters on opposite sides of the threshold. A further extension is to treat the
initial densities as unknown and to study how uncertainty in the initial
infected profile affects inference of persistence.
\par
\endgroup

A related question is what can be inferred from case reports collected over time intervals. In an earlier work on SIR \cite{CampilloFunolletEtAl2022}, we treated reports as a fraction of the exits from the infected class.
For a spatial SIS analogue, let $r(x)$ denote the fraction of local transitions from $I$ to $S$ that are reported. The cumulative density of reported cases $D$ then satisfies $J=\partial_tD=\kappa I$, with $\kappa=r\gamma$. If $\kappa$ is bounded away from zero, we can reformulate \eqref{eq:SIS} in terms of $N$ and $J$ by substituting $I=J/\kappa$ and $S=N-J/\kappa$. The diffusion term and Neumann condition for the infected density then involve $J/\kappa$, and $N$ and $J$ remain coupled when $d_S\ne d_I$. Integrating $D(t)-D(s)$ over an observation region gives the modelled count over $(s,t]$. Extending the growth bounds to these counts would help relate interval observations to the persistence threshold.
This raises the question of how to choose the observation regions and reporting intervals to minimise expected posterior uncertainty about whether \(\mathcal R_0>1\). This is a classic problem in \emph{Bayesian optimal experimental design}, which we plan to study in future work.

\section*{Acknowledgment}
The author would like to thank Tapio Helin for suggesting the study of PDE models of infectious diseases and sharing
\cite{CoronelFrizHessZegarra2021}, which led to the development of this project. This work is partly supported by the Research Council of Finland via FAME Flagship (decision numbers 359183 and 359186) and the Magnus Ehrnrooth Foundation. Part of this research was performed while the author was visiting the Institute for Mathematical and Statistical Innovation (IMSI), which is supported by the National Science Foundation (Grant No. DMS-2425650). The author would like to thank the organizers of the IMSI Fall 2025 program \emph{Digital Twins: Mathematical and Statistical Foundations and Complex Applications}, during which a substantial part of this work was done.

\section*{AI disclosure} During the preparation of this manuscript, ChatGPT (OpenAI) was used to assist with literature exploration, proof-reading, presentation and to help implement the numerical experiments. All the AI-assisted outputs were reviewed and independently verified by the author, who takes full responsibility for the content of this paper.

\appendix
\crefalias{section}{appendix}
\crefalias{subsection}{appendix}
\section{Uniform spectral constants and weighted observations}\label{app:spectral}
We first derive the uniform constants in \cref{prop:uniform-spectral}.
Fix $\tau_0>0$ and let $K_d(\tau_0,x,y)$ be the kernel of the Neumann
heat semigroup $P_d(\tau_0)$ in \cref{sec:heat-evolution}. The kernel
is continuous on the closed domain at positive time
\cite[Sections~3.1 and~4.2]{ArendtTerElstGluck2020}. For each $y$,
$K_d(\tau_0/2,\cdot,y)$ is nonnegative and has integral one by symmetry
and preservation of constants. The semigroup property and positivity
improvement then give $K_d(\tau_0,x,y)>0$ on
$\overline\Omega\times\overline\Omega$. The time scaling
$K_d(\tau_0,x,y)=K_1(d\tau_0,x,y)$ and continuity at positive times give
\[
 k_-:=\min_{\substack{x,y\in\overline\Omega\\d\in[\underline d,\overline d]}}
 K_d(\tau_0,x,y)>0,\qquad
 k_+:=\max_{\substack{x,y\in\overline\Omega\\d\in[\underline d,\overline d]}}
 K_d(\tau_0,x,y)<\infty.
\]
For any nonnegative $f$, comparison with the constant potentials $q_-,q_+$
gives
\[
 e^{q_-\tau_0}P_{d_I}(\tau_0)f
 \le e^{\tau_0 L_\theta}f
 \le e^{q_+\tau_0}P_{d_I}(\tau_0)f.
\]
We apply these inequalities to $f=\phi_1$, for which
$e^{\tau_0L_\theta}\phi_1=e^{\lambda_1\tau_0}\phi_1$.
Taking the upper bound at one point and the lower bound at another gives
\begin{equation}\label{eq:kernel-H}
 \frac{\max\phi_1}{\min\phi_1}
 \le \mathfrak h:=e^{(q_+-q_-)\tau_0}\frac{k_+}{k_-}.
\end{equation}
The constant $\mathfrak h$ depends only on the coefficient bounds,
the domain, and the chosen comparison time.
Normalization in $L^2$ implies
$\min\phi_1\ge1/(\mathfrak h\sqrt V)$ and
$\max\phi_1\le \mathfrak h/\sqrt V$, proving the overlap bound in
\cref{prop:uniform-spectral}.

For completeness, integration by parts in the principal eigenfunction equation
gives the ground-state identity
\[
 \lambda_1\|\phi_1v\|_{L^2}^2
 -\langle L_\theta(\phi_1v),\phi_1v\rangle
 =d_I\int_\Omega\phi_1^2|\nabla v|^2\dd x.
\]
By density and the variational characterization of $\lambda_2$, the gap
is the infimum of the associated quotient over
$\int\phi_1^2v=0$. For such $v$, since $\int\phi_1^2=1$,
\begin{align*}
 \int\phi_1^2v^2
 &=\inf_{c\in\R}\int\phi_1^2(v-c)^2
 \le\frac{\max\phi_1^2}{\nu_1}\int|\nabla v|^2\\
 &\le\frac{\mathfrak h^2}{\nu_1}\int\phi_1^2|\nabla v|^2.
\end{align*}
It follows that $\mathfrak g_\theta\ge d_I\nu_1/\mathfrak h^2$, which
completes the proof of the uniform spectral estimates.

On $\Omega=(0,L)$, the cosine heat-kernel expansion gives usable constants.
Choose $\tau_0$ so that
$z=2\sum_{n\ge1}e^{-\underline d(n\pi/L)^2\tau_0}<1$.
Then $k_-\ge(1-z)/L$, $k_+\le(1+z)/L$, and
$\nu_1=\pi^2/L^2$. Thus
\[
 \mathfrak h\le e^{(q_+-q_-)\tau_0}\frac{1+z}{1-z}.
\]
The positive series has elementary integral tail bounds. These constants
can be conservative. Sharper verified heat-kernel or eigenfunction ratio
bounds can be substituted in the estimates above.

\begin{proposition}\label{prop:weighted-growth}
Under the hypotheses of \cref{thm:growth-threshold}, let
$0\ne w\in L^2(\Omega)$ be nonnegative, let
$b_1=\langle w,\phi_1\rangle>0$, and set
\[
 D_w(u)=\frac{\rho_\theta\|w\|_{L^2}}{b_1}
          e^{-\mathfrak g_\theta u},\qquad
 \Lambda_{w,\theta}(s,t)=\frac1{t-s}
 \log\frac{\langle i_\theta(t),w\rangle}{\langle i_\theta(s),w\rangle}.
\]
If $D_w(s)<1$, then
\begin{equation}\label{eq:weighted-growth-bound}
 |\Lambda_{w,\theta}(s,t)-\lambda_*|
 \le\frac{-\log(1-D_w(s))-\log(1-D_w(t))}{t-s}=:W_{w,\theta}(s,t).
\end{equation}
For the nonlinear statistic $\Lambda_{w,\theta}^\varepsilon$ formed from
$I_\theta^\varepsilon$, an envelope as in \cref{thm:growth-threshold} gives
\[
 |\Lambda_{w,\theta}^\varepsilon(s,t)-\lambda_*|
 \le W_{w,\theta}(s,t)+B_U(s,t).
\]
The constants can be made uniform: use
$b_1\ge\|w\|_{L^1}/(\mathfrak h\sqrt V)$ and \eqref{eq:uniform-spectral}.
\end{proposition}
\begin{proof}
By the eigenfunction expansion,
\[
 \langle i_\theta(u),w\rangle=c_1b_1e^{\lambda_1u}(1+r_w(u)).
\]
By the Cauchy--Schwarz inequality and the spectral gap,
$|r_w(u)|\le D_w(u)$. For $|r|\le D<1$, it follows that
$|\log(1+r)|\le-\log(1-D)$. Taking the logarithmic ratio, we obtain
\eqref{eq:weighted-growth-bound}.
By \eqref{eq:multiplicative-comparison}, the weighted nonlinear-to-linear
log ratio at time $u$ belongs to $[-A_U(u),0]$. Its difference between
$s$ and $t$ has absolute value at most $A_U(t)$, proving the nonlinear
claim. The lower bound for $b_1$ follows from \eqref{eq:kernel-H}.
\end{proof}

For $w=1$, these are total infected counts. Unlike the squared-norm expansion
\eqref{eq:growth-spectral-expansion}, the weighted spectral remainder can
have either sign. A persistence conclusion based on weighted counts
therefore requires the corresponding transient-error allowance.
The same logarithmic measurement-error term $d_{\rm obs}$ can be added
when weighted counts are observed with the stated relative error bound.

\section{Prior construction and approximation}\label{app:prior-approximation}

\subsection{Latent prior and truncation}\label{app:latent-truncation}
Let $\nu_0=N(m,\mathcal C)$ on $\mathcal U$ and let $\mathcal T$ be an
admissible Lipschitz transformation from \cref{sec:prior-modelling}. The
latent posterior is given by
\begin{equation}\label{eq:latent-posterior}
 \frac{\dd\nu^y}{\dd\nu_0}(u)=Z(y)^{-1}e^{-\Psi(u;y)},
 \qquad \Psi(u;y)=\Phi(\mathcal T(u);y).
\end{equation}
Integrating a bounded Borel function $f(\mathcal T(u))$ against this
density gives the same numerator and normalizing constant as integrating
$f(\theta)$ against \eqref{eq:posterior-RN}, since
$\mu_0=\mathcal T_\sharp\nu_0$. Thus $\mathcal T_\sharp\nu^y=\mu^y$.

On $\Omega=(0,L)$, let $e_0=L^{-1/2}$ and
$e_j=(2/L)^{1/2}\cos(j\pi x/L)$ for $j\ge1$. A continuous latent field
can be constructed from the series
\begin{equation}\label{eq:cosine-prior}
 u(x)=m(x)+\sum_{j=0}^\infty a_j\xi_j e_j(x),\qquad
 a_j=\sigma\bigl(\kappa^2+(j\pi/L)^2\bigr)^{-p/2},\qquad
 \xi_j\overset{\mathrm{iid}}{\sim}N(0,1),
\end{equation}
where $m\in C([0,L])$, $\sigma,\kappa>0$, and $p>1$. The pointwise
variance is $\sum_{j\ge0}a_j^2e_j(x)^2$. Since the basis is uniformly
bounded and $\sum_j a_j<\infty$, the expected sum of the sup norms of
the terms is finite. The series therefore converges absolutely and
uniformly almost surely. For the truncation $u_N$ at mode $N$,
Minkowski's inequality gives
\begin{equation}\label{eq:latent-tail}
 \bigl(\E\norm{u-u_N}_\infty^2\bigr)^{1/2}
 \le\sum_{j>N}a_j\norm{e_j}_\infty\le C N^{1-p}.
\end{equation}
Two such fields can be generated independently or coupled by a fixed
linear mixing before applying $\mathcal T$.

For a measurable approximation $u_N$ of $u\sim\nu_0$, suppose that
$\E\norm{u-u_N}_{\mathcal U}^2\le\epsilon_N^2$.
Define $\nu_N^y$ by replacing $\Psi(u;y)$ with $\Psi(u_N;y)$ in
\eqref{eq:latent-posterior} and normalizing against the same $\nu_0$.
The Lipschitz property of $\mathcal T$ and the estimate in
\cref{prop:G-regularity} give an $L^2(\nu_0)$ observation error of
order $\epsilon_N$. The proof of \cref{thm:posterior-stability}
uses only a common prior probability space and bounded measurable
observation maps. Applying it on $(\mathcal U,\nu_0)$ gives
\begin{equation}\label{eq:latent-truncation-Hell}
 \sup_{|y|\le R}d_\Hell(\nu^y,\nu_N^y)\le C_R\epsilon_N.
\end{equation}
This estimate compares densities on the common latent space. For the series
truncation in \eqref{eq:cosine-prior}, the likelihood depends only on the
retained coordinates, so the independent tail keeps its prior law under
$\nu_N^y$. The retained coordinates can be sampled using preconditioned
Crank--Nicolson with potential $\Psi(u_N;y)$
\cite{CotterRobertsStuartWhite2013}.

\subsection{Forward approximation}\label{app:forward}
Let $P_N:\X_\ad\to\X_\ad$ be measurable and satisfy
$P_N\theta\to\theta$ in $\X$ for $\mu_0$-almost every $\theta$. Define
\begin{equation}\label{eq:typical-forward-rate}
 \epsilon_N=\left(\int_{\X_\ad}
       \norm{\theta-P_N\theta}_{\X}^2\dd\mu_0(\theta)\right)^{1/2}.
\end{equation}
Boundedness of $\X_\ad$ gives $\epsilon_N\to0$ by dominated convergence.
By \cref{prop:G-regularity}, the error from replacing $\G$ by
$\G\circ P_N$ in $L^2(\mu_0;\R^J)$ is at most $C\epsilon_N$.
Suppose a solver with mesh size $h$ and time step $\tau$ gives a measurable,
uniformly bounded approximate observation map. If \eqref{eq:state-solver-error}
at the projected coefficients holds with error $C_N(\theta)(h^r+\tau)$
and $\sup_N\norm{C_N}_{L^2(\mu_0)}<\infty$, then
\cref{thm:posterior-approximation} gives
\begin{equation}\label{eq:typical-posterior-rate}
 \sup_{|y|\le R}d_\Hell(\mu^y,\mu_{N,h,\tau}^y)
 \le C_R(\epsilon_N+h^r+\tau).
\end{equation}
The solver estimate requires stability and spatial and temporal
regularity appropriate to the discretization. For parabolic finite-element
estimates, see \cite{Thomee2006}. The next result supplies an eigenvalue
error bound to use together with \eqref{eq:typical-posterior-rate}
in \cref{thm:threshold-approx}.

\subsection{Conforming approximation of the threshold}
For a conforming finite-element space $V_h\subset H^1(\Omega)$, define
$\Rzero{}_h$ by restricting the supremum in \eqref{eq:R0-def} to $V_h$.
Let $\phi_R(\theta)$ be the positive principal eigenfunction, normalized by
$d_I\|\nabla\phi_R\|_{L^2}^2+\int_\Omega\gamma\phi_R^2\dd x=1$.

\begin{corollary}\label{cor:eigen-FE}
Suppose there are $r>0$ and $C<\infty$ such that
\begin{equation}\label{eq:eigen-approx-assump}
 \inf_{v_h\in V_h}\norm{\phi_R(\theta)-v_h}_{H^1}
 \le Ch^r\qquad\text{for all }\theta\in\X_\ad.
\end{equation}
Then the conforming approximation satisfies
\begin{equation}\label{eq:R0-FE-rate}
 0\le\Rzero(\theta)-\Rzero{}_h(\theta)\le Ch^{2r}
 \qquad\text{uniformly over }\X_\ad.
\end{equation}
If the forward posterior error is $O(\psi_{\mathfrak a})$ in Hellinger
distance, the error in the posterior mean of $\Rzero$ is
$O(\psi_{\mathfrak a}+h^{2r})$. The approximate persistence probability
converges under the no-threshold-atom condition of \cref{thm:threshold-approx}.
\end{corollary}
\begin{proof}
Write $a_\theta$ and $b_\theta$ for the denominator and numerator forms
of \eqref{eq:R0-def}, respectively. Conformity gives $\Rzero{}_h\le\Rzero$.
For $v_h$ with $\|v_h-\phi_R\|_{H^1}\le Ch^r$, set $e_h=v_h-\phi_R$.
The identity $b_\theta(\phi_R,v)=\Rzero\,a_\theta(\phi_R,v)$ cancels the
linear terms in the Rayleigh-quotient error \cite{Boffi2010}, giving
\[
 \Rzero-\frac{b_\theta(v_h,v_h)}{a_\theta(v_h,v_h)}
 =\frac{\Rzero\,a_\theta(e_h,e_h)-b_\theta(e_h,e_h)}
        {a_\theta(v_h,v_h)}.
\]
The numerator is $O(h^{2r})$ by the uniform bounds on the forms and
$\Rzero$. Since $a_\theta(\phi_R,\phi_R)=1$, the denominator is bounded
away from zero for sufficiently small $h$, uniformly in $\theta$.
Taking the discrete supremum gives \eqref{eq:R0-FE-rate}.
The posterior conclusions follow from \cref{thm:threshold-approx}.
\end{proof}

On the fixed $C^2$ domain the normalized eigenfunction has a uniform
$H^2$ bound. Indeed, $-d_I\Delta\phi_R=(\beta/\Rzero-\gamma)\phi_R$ with
Neumann boundary conditions, and
$\underline\beta/\overline\gamma\le\Rzero\le\overline\beta/\underline\gamma$.
Thus $r=1$ applies to conforming spaces with first-order $H^1$
approximation of $H^2$ functions, with boundary geometry errors treated
separately. The finite-volume schemes used for the figures are described
in \cref{app:implementation}.

\section{Details of numerical experiments}\label{app:implementation}
We use the parameterization in \eqref{eq:numerical-parameterization} and
independent Gaussian priors truncated to the indicated intervals
\begin{align*}
 q_0&\sim N(-0.015,0.025^2),& q_0&\in[-0.060,0.030],\\*
 q_1&\sim N(0.080,0.045^2),& q_1&\in[0,0.180],\\*
 a&\sim N(0.360,0.050^2),& a&\in[0.280,0.460].
\end{align*}
We use a cell-centred finite-volume discretization on uniform cells of
width $h$. Initial profiles and coefficients are evaluated at the
centres. Zero boundary flux gives the tridiagonal Laplacian $D_h$
with off-diagonal entries $h^{-2}$ and diagonal entries $-2h^{-2}$ in
the interior and $-h^{-2}$ at the two boundary cells. Observations are
averages on cells with one-based indices
$(3,6,10,13,17,21,24,28,31,35,39)$ of the 41-cell grid, at $t=0.5$
and $t=2$. We use the noise standard deviations given in
\cref{sec:numerics}. The normal draws use NumPy's
\texttt{default\_rng}, ordered by time then cell, with prevalence
followed by incidence.

\subsection{Controlled limiting-model experiment}
For \cref{fig:posterior-results}, both data generation and inference
use the 41-cell operator $L_h=d_ID_h+\operatorname{diag}(q)$.
We evaluate its matrix exponential by symmetric eigendecomposition.
The observation means are the selected components of $i_h(t)$ and
$\beta i_h(t)$, corresponding to \cref{cor:joint-limit}.
We generate the noise with seed $24051983$. Posterior masses are
proportional to the products of the prior densities, the Gaussian
likelihood, and tensor trapezoidal weights on a
$121\times121\times101$ grid in $(q_0,q_1,a)$. The recovery intervals
use the cumulative trapezoidal integral of the marginal density.
To obtain the eigenvalue density, we use the shift identity
$\lambda_{*,h}(q_0,q_1)=q_0+\lambda_{*,h}(0,q_1)$. For each plotted
eigenvalue, we interpolate the $(q_0,q_1)$ marginal density in $q_0$
at the corresponding shifted value and integrate in $q_1$. The threshold
probabilities and the growth bracket in \cref{sec:numerics} use
quadrature of their event indicators, with the discrete $L^2$ norm
and spectral decomposition of $L_h$ for the growth calculation.

\subsection{Nonlinear posterior comparison}\label{app:nonlinear-posterior}
For \cref{fig:posterior-limit}, data are generated at the true parameter
on 123 cells by a backward differentiation formula (BDF) in
$(S,I/\varepsilon)$, with relative and absolute tolerances
$2\times10^{-10}$ and $2\times10^{-12}$. The means are averages of
$I/\varepsilon$ and $\beta\F(S,I)/\varepsilon$ over the three fine cells
within each observation cell. For each $\varepsilon$, we add the same
standard normal vector generated with seed $20261104$, scaled by the
observation standard deviations.

Both likelihoods use the 41-cell grid. For the nonlinear solve we use
Strang splitting with step $\Delta t=0.02$, applying each diffusion
half-step by the exact discrete Neumann heat semigroup. In the reaction
step $N=S+I$ is constant on each cell, and the update is
\[
 I_{\rm new}=\frac{I e^{q\Delta t}}
 {1+(\beta I/N)\Delta t\,\operatorname{exprel}(q\Delta t)},\qquad
 \operatorname{exprel}(z)=\frac{e^z-1}{z},\quad
 \operatorname{exprel}(0)=1 .
\]
We then set $S_{\rm new}=N-I_{\rm new}$. The limiting solve uses the
matrix exponential as above. With normalized posterior masses $p_k,q_k$
on a common $41\times41\times31$ trapezoidal grid, the plotted distance is
\begin{equation}\label{eq:discrete-hellinger}
 H_{\rm grid}=\left[\frac12\sum_k(\sqrt{p_k}-\sqrt{q_k})^2\right]^{1/2}.
\end{equation}
At $\varepsilon=0.00125$, the distances are $0.030401$ for prevalence
and $0.057887$ for joint data. Changing the parameter grid to
$31\times31\times21$ or $61\times61\times41$, halving $\Delta t$, or
doubling the inversion mesh to 82 cells changes either distance by
at most $2.3\times10^{-6}$ in these checks. On the 82-cell grid the
observations average the two cells in each original observation interval.

For the state errors in \cref{fig:posterior-limit}, we use the true
parameter and the 41-cell grid for all states. The nonlinear reference
uses the BDF tolerances above. The correction in \eqref{eq:second-order}
has forcing $-\beta i_h^2$ since $S^{(0)}=1$. We integrate it with BDF
using relative and absolute tolerances $2\times10^{-11}$ and $2\times10^{-13}$.
The errors at $t=2$ are $\|I^\varepsilon/\varepsilon-i_h\|_h$ and
$\|I^\varepsilon/\varepsilon-i_h-\varepsilon j_h\|_h$, where
$\|v\|_h^2=h\sum_j|v_j|^2$. Log--log slope fits use the three smallest
infection scales for the Hellinger distances and all seven for the
state errors.

\bibliographystyle{siamplain}
\bibliography{BIPs_for_SIS-arXiv-v3-4}

@article{AllenBolkerLouNevai2008,
  author  = {Allen, Linda J. S. and Bolker, Benjamin M. and Lou, Yuan and Nevai, Alexander L.},
  title   = {Asymptotic Profiles of the Steady States for an {SIS} Epidemic Reaction--Diffusion Model},
  journal = {Discrete Contin. Dyn. Syst.},
  volume  = {21},
  number  = {1},
  pages   = {1--20},
  year    = {2008},
  doi     = {10.3934/dcds.2008.21.1}
}

@article{ArendtTerElstGluck2020,
  author  = {Arendt, Wolfgang and ter Elst, A. F. M. and Gl{\"u}ck, Jochen},
  title   = {Strict positivity for the principal eigenfunction of elliptic operators with various boundary conditions},
  journal = {Adv. Nonlinear Stud.},
  volume  = {20},
  number  = {3},
  pages   = {633--650},
  year    = {2020},
  doi     = {10.1515/ans-2020-2091}
}

@article{CampilloFunolletEtAl2022,
  author  = {Campillo-Funollet, Eduard and Wragg, Hayley and Van Yperen, James and Duong, Duc-Lam and Madzvamuse, Anotida},
  title   = {Reformulating the susceptible--infectious--removed model in terms of the number of detected cases: well-posedness of the observational model},
  journal = {Philos. Trans. R. Soc. A},
  volume  = {380},
  number  = {2233},
  pages   = {20210306},
  year    = {2022},
  doi     = {10.1098/rsta.2021.0306}
}

@article{Peng2009,
  author  = {Peng, Rui},
  title   = {Asymptotic Profiles of the Positive Steady State for an {SIS} Epidemic Reaction--Diffusion Model. {Part I}},
  journal = {J. Differential Equations},
  volume  = {247},
  number  = {4},
  pages   = {1096--1119},
  year    = {2009},
  doi     = {10.1016/j.jde.2009.05.002}
}

@article{PengYi2013,
  author  = {Peng, Rui and Yi, Fengqi},
  title   = {Asymptotic profile of the positive steady state for an {SIS} epidemic reaction--diffusion model: Effects of epidemic risk and population movement},
  journal = {Phys. D},
  volume  = {259},
  pages   = {8--25},
  year    = {2013},
  doi     = {10.1016/j.physd.2013.05.006}
}

@article{CuiLou2016,
  author  = {Cui, Renhao and Lou, Yuan},
  title   = {A spatial {SIS} model in advective heterogeneous environments},
  journal = {J. Differential Equations},
  volume  = {261},
  number  = {6},
  pages   = {3305--3343},
  year    = {2016},
  doi     = {10.1016/j.jde.2016.05.025}
}

@article{LiPengWang2017,
  author  = {Li, Huicong and Peng, Rui and Wang, Feng-Bin},
  title   = {Varying total population enhances disease persistence: Qualitative analysis on a diffusive {SIS} epidemic model},
  journal = {J. Differential Equations},
  volume  = {262},
  number  = {2},
  pages   = {885--913},
  year    = {2017},
  doi     = {10.1016/j.jde.2016.09.044}
}

@article{HuJinWang2024,
  author  = {Hu, Yaru and Jin, Yu and Wang, Jinfeng},
  title   = {Dynamics of a Reaction-Diffusion {SIS} Epidemic Model with a Control Zone},
  journal = {SIAM J. Appl. Math.},
  volume  = {84},
  number  = {6},
  pages   = {2569--2589},
  year    = {2024},
  doi     = {10.1137/24M1652374}
}

@article{HuangHanLiu2010,
  author  = {Huang, Wenzhang and Han, Maoan and Liu, Kaiyu},
  title   = {Dynamics of an {SIS} Reaction--Diffusion Epidemic Model for Disease Transmission},
  journal = {Math. Biosci. Eng.},
  volume  = {7},
  number  = {1},
  pages   = {51--66},
  year    = {2010},
  doi     = {10.3934/mbe.2010.7.51}
}

@article{MagalWebbWu2019,
  author  = {Magal, Pierre and Webb, Glenn F. and Wu, Yixiang},
  title   = {On the Basic Reproduction Number of Reaction--Diffusion Epidemic Models},
  journal = {SIAM J. Appl. Math.},
  volume  = {79},
  number  = {1},
  pages   = {284--304},
  year    = {2019},
  doi     = {10.1137/18M1182243}
}

@article{XiangLiu2015,
  author  = {Xiang, Huili and Liu, Bin},
  title   = {Solving the Inverse Problem of an {SIS} Epidemic Reaction--Diffusion Model by Optimal Control Methods},
  journal = {Comput. Math. Appl.},
  volume  = {70},
  number  = {5},
  pages   = {805--819},
  year    = {2015},
  doi     = {10.1016/j.camwa.2015.05.025}
}

@article{CoronelFrizHessZegarra2021,
  author  = {Coronel, An{\'i}bal and Friz, Luis and Hess, Ian and Zegarra, Mar{\'i}a},
  title   = {On the Existence and Uniqueness of an Inverse Problem in Epidemiology},
  journal = {Appl. Anal.},
  volume  = {100},
  number  = {3},
  pages   = {513--526},
  year    = {2021},
  doi     = {10.1080/00036811.2019.1608964}
}

@article{CoronelHuancasSepulveda2019,
  author  = {Coronel, An{\'i}bal and Huancas, Fernando and Sep{\'u}lveda, Mauricio},
  title   = {A Note on the Existence and Stability of an Inverse Problem for a {SIS} Model},
  journal = {Comput. Math. Appl.},
  volume  = {77},
  number  = {12},
  pages   = {3186--3194},
  year    = {2019},
  doi     = {10.1016/j.camwa.2019.01.031}
}

@article{CoronelHuancasHessTello2024,
  author  = {Coronel, An{\'i}bal and Huancas, Fernando and Hess, Ian and Tello, Alex},
  title   = {The Diffusion Identification in a {SIS} Reaction--Diffusion System},
  journal = {Math. Biosci. Eng.},
  volume  = {21},
  number  = {1},
  pages   = {562--581},
  year    = {2024},
  doi     = {10.3934/mbe.2024024}
}

@article{JhaCaoOden2020,
  author  = {Jha, Prashant K. and Cao, Lianghao and Oden, J. Tinsley},
  title   = {Bayesian-Based Predictions of {COVID-19} Evolution in {Texas} Using Multispecies Mixture-Theoretic Continuum Models},
  journal = {Comput. Mech.},
  volume  = {66},
  number  = {5},
  pages   = {1055--1068},
  year    = {2020},
  doi     = {10.1007/s00466-020-01889-z}
}

@article{ChenWuGhattas2021,
  author  = {Chen, Peng and Wu, Keyi and Ghattas, Omar},
  title   = {Bayesian Inference of Heterogeneous Epidemic Models: Application to {COVID-19} Spread Accounting for Long-Term Care Facilities},
  journal = {Comput. Methods Appl. Mech. Engrg.},
  volume  = {385},
  pages   = {114020},
  year    = {2021},
  doi     = {10.1016/j.cma.2021.114020}
}

@article{CotterDashtiRobinsonStuart2009,
  author  = {Cotter, Simon L. and Dashti, Masoumeh and Robinson, James C. and Stuart, Andrew M.},
  title   = {Bayesian Inverse Problems for Functions and Applications to Fluid Mechanics},
  journal = {Inverse Problems},
  volume  = {25},
  number  = {11},
  pages   = {115008},
  year    = {2009},
  doi     = {10.1088/0266-5611/25/11/115008}
}

@article{Stuart2010,
  author  = {Stuart, Andrew M.},
  title   = {Inverse Problems: A {Bayesian} Perspective},
  journal = {Acta Numer.},
  volume  = {19},
  pages   = {451--559},
  year    = {2010},
  doi     = {10.1017/S0962492910000061}
}

@article{CotterDashtiStuart2010,
  author  = {Cotter, Simon L. and Dashti, Masoumeh and Stuart, Andrew M.},
  title   = {Approximation of {Bayesian} Inverse Problems for {PDEs}},
  journal = {SIAM J. Numer. Anal.},
  volume  = {48},
  number  = {1},
  pages   = {322--345},
  year    = {2010},
  doi     = {10.1137/090770734}
}

@article{Latz2020,
  author  = {Latz, Jonas},
  title   = {On the Well-Posedness of {Bayesian} Inverse Problems},
  journal = {SIAM/ASA J. Uncertain. Quantif.},
  volume  = {8},
  number  = {1},
  pages   = {451--482},
  year    = {2020},
  doi     = {10.1137/19M1247176}
}

@incollection{DashtiStuart2017,
  author    = {Dashti, Masoumeh and Stuart, Andrew M.},
  title     = {The {Bayesian} Approach to Inverse Problems},
  booktitle = {Handbook of Uncertainty Quantification},
  editor    = {Ghanem, Roger and Higdon, David and Owhadi, Houman},
  pages     = {311--428},
  publisher = {Springer International Publishing},
  address   = {Cham},
  year      = {2017},
  doi       = {10.1007/978-3-319-12385-1_7}
}

@book{Lunardi1995,
  author    = {Lunardi, Alessandra},
  title     = {Analytic Semigroups and Optimal Regularity in Parabolic Problems},
  series    = {Progress in Nonlinear Differential Equations and Their Applications},
  volume    = {16},
  publisher = {Birkh{\"a}user},
  address   = {Basel},
  year      = {1995}
}

@book{Thomee2006,
  author    = {Thom{\'e}e, Vidar},
  title     = {Galerkin Finite Element Methods for Parabolic Problems},
  edition   = {Second},
  series    = {Springer Series in Computational Mathematics},
  volume    = {25},
  publisher = {Springer},
  address   = {Berlin},
  year      = {2006},
  doi       = {10.1007/3-540-33122-0}
}

@article{Boffi2010,
  author  = {Boffi, Daniele},
  title   = {Finite Element Approximation of Eigenvalue Problems},
  journal = {Acta Numer.},
  volume  = {19},
  pages   = {1--120},
  year    = {2010},
  doi     = {10.1017/S0962492910000012}
}

@article{LindgrenRueLindstrom2011,
 author = {Lindgren, Finn and Rue, H{\aa}vard and Lindstr{\"o}m, Johan},
 title = {An Explicit Link between {Gaussian} Fields and {Gaussian} {Markov} Random Fields: The Stochastic Partial Differential Equation Approach},
 journal = {J. R. Stat. Soc. Ser. B. Stat. Methodol.},
 volume = {73}, number = {4}, pages = {423--498}, year = {2011},
 doi = {10.1111/j.1467-9868.2011.00777.x}
}

@article{CotterRobertsStuartWhite2013,
 author = {Cotter, Simon L. and Roberts, Gareth O. and Stuart, Andrew M. and White, David},
 title = {{MCMC} Methods for Functions: Modifying Old Algorithms to Make Them Faster},
 journal = {Statist. Sci.}, volume = {28}, number = {3},
 pages = {424--446}, year = {2013}, doi = {10.1214/13-STS421}
}

@article{Poirier1998,
 author={Poirier, Dale J.},
 title={Revising Beliefs in Nonidentified Models},
 journal={Econometric Theory},volume={14},number={4},pages={483--509},year={1998},
 doi={10.1017/S0266466698144043}
}

@book{Gustafson2015,
 author={Gustafson, Paul},
 title={Bayesian Inference for Partially Identified Models: Exploring the Limits of Limited Data},
 publisher={Chapman and Hall/CRC},year={2015},doi={10.1201/b18308}
}

@article{WangZhao2012,
 author={Wang, Wendi and Zhao, Xiao-Qiang},
 title={Basic Reproduction Numbers for Reaction--Diffusion Epidemic Models},
 journal={SIAM J. Appl. Dyn. Syst.},volume={11},number={4},pages={1652--1673},
 year={2012},doi={10.1137/120872942}
}

@article{WallingaLipsitch2007,
 author={Wallinga, Jacco and Lipsitch, Marc},
 title={How Generation Intervals Shape the Relationship between Growth Rates and Reproductive Numbers},
 journal={Proc. R. Soc. B},volume={274},number={1609},pages={599--604},
 year={2007},doi={10.1098/rspb.2006.3754}
}

@article{Sprungk2020,
 author={Sprungk, Bj{\"o}rn},
 title={On the Local {Lipschitz} Stability of {Bayesian} Inverse Problems},
 journal={Inverse Problems},volume={36},number={5},pages={055015},
 year={2020},doi={10.1088/1361-6420/ab6f43}
}

@article{hethcote2000mathematics,
  title={The mathematics of infectious diseases},
  author={Hethcote, Herbert W},
  journal={SIAM review},
  volume={42},
  number={4},
  pages={599--653},
  year={2000},
  publisher={SIAM}
}
\end{document}